\documentclass[12pt]{article}
\usepackage{latexsym, amssymb}
\usepackage{dsfont}

\DeclareMathAlphabet\gothic{U}{euf}{m}{n}

\makeatletter
\def\eqnarray{\stepcounter{equation}\let\@currentlabel=\theequation
\global\@eqnswtrue
\tabskip\@centering\let\\=\@eqncr
$$\halign to \displaywidth\bgroup\hfil\global\@eqcnt\z@
  $\displaystyle\tabskip\z@{##}$&\global\@eqcnt\@ne
  \hfil$\displaystyle{{}##{}}$\hfil
  &\global\@eqcnt\tw@ $\displaystyle{##}$\hfil
  \tabskip\@centering&\llap{##}\tabskip\z@\cr}

\def\endeqnarray{\@@eqncr\egroup
      \global\advance\c@equation\m@ne$$\global\@ignoretrue}

\def\@yeqncr{\@ifnextchar [{\@xeqncr}{\@xeqncr[5pt]}}
\makeatother

\begin{document}
\bibliographystyle{tom}

\newtheorem{lemma}{Lemma}[section]
\newtheorem{thm}[lemma]{Theorem}
\newtheorem{cor}[lemma]{Corollary}
\newtheorem{voorb}[lemma]{Example}
\newtheorem{rem}[lemma]{Remark}
\newtheorem{prop}[lemma]{Proposition}
\newtheorem{ddefinition}[lemma]{Definition}
\newtheorem{stat}[lemma]{{\hspace{-5pt}}}

\newenvironment{remarkn}{\begin{rem} \rm}{\end{rem}}
\newenvironment{exam}{\begin{voorb} \rm}{\end{voorb}}
\newenvironment{definition}{\begin{ddefinition} \rm}{\end{ddefinition}}

\newcommand{\gota}{\gothic{a}}
\newcommand{\gotb}{\gothic{b}}
\newcommand{\gotc}{\gothic{c}}
\newcommand{\gote}{\gothic{e}}
\newcommand{\gotf}{\gothic{f}}
\newcommand{\gotg}{\gothic{g}}
\newcommand{\gothh}{\gothic{h}}
\newcommand{\gotk}{\gothic{k}}
\newcommand{\gotm}{\gothic{m}}
\newcommand{\gotn}{\gothic{n}}
\newcommand{\gotp}{\gothic{p}}
\newcommand{\gotq}{\gothic{q}}
\newcommand{\gotr}{\gothic{r}}
\newcommand{\gots}{\gothic{s}}
\newcommand{\gotu}{\gothic{u}}
\newcommand{\gotv}{\gothic{v}}
\newcommand{\gotw}{\gothic{w}}
\newcommand{\gotz}{\gothic{z}}
\newcommand{\gotA}{\gothic{A}}
\newcommand{\gotB}{\gothic{B}}
\newcommand{\gotG}{\gothic{G}}
\newcommand{\gotL}{\gothic{L}}
\newcommand{\gotS}{\gothic{S}}
\newcommand{\gotT}{\gothic{T}}

\newcounter{teller}
\renewcommand{\theteller}{(\alph{teller})}
\newenvironment{tabel}{\begin{list}%
{\rm  (\alph{teller})\hfill}{\usecounter{teller} \leftmargin=1.1cm
\labelwidth=1.1cm \labelsep=0cm \parsep=0cm}
                      }{\end{list}}

\newcounter{tellerr}
\renewcommand{\thetellerr}{(\roman{tellerr})}
\newenvironment{tabeleq}{\begin{list}%
{\rm  (\roman{tellerr})\hfill}{\usecounter{tellerr} \leftmargin=1.1cm
\labelwidth=1.1cm \labelsep=0cm \parsep=0cm}
                         }{\end{list}}

\newcounter{tellerrr}
\renewcommand{\thetellerrr}{(\Roman{tellerrr})}
\newenvironment{tabelR}{\begin{list}%
{\rm  (\Roman{tellerrr})\hfill}{\usecounter{tellerrr} \leftmargin=1.1cm
\labelwidth=1.1cm \labelsep=0cm \parsep=0cm}
                         }{\end{list}}

\newcounter{proofstep}
\newcommand{\nextstep}{\refstepcounter{proofstep}\vertspace \par 
          \noindent{\bf Step \theproofstep} \hspace{5pt}}
\newcommand{\firststep}{\setcounter{proofstep}{0}\nextstep}

\newcommand{\Ni}{\mathds{N}}
\newcommand{\Qi}{\mathds{Q}}
\newcommand{\Ri}{\mathds{R}}
\newcommand{\Ci}{\mathds{C}}
\newcommand{\Ti}{\mathds{T}}
\newcommand{\Zi}{\mathds{Z}}
\newcommand{\Fi}{\mathds{F}}

\newcommand{\proof}{\mbox{\bf Proof} \hspace{5pt}} 
\newcommand{\remark}{\mbox{\bf Remark} \hspace{5pt}}
\newcommand{\vertspace}{\vskip10.0pt plus 4.0pt minus 6.0pt}

\newcommand{\simh}{{\stackrel{{\rm cap}}{\sim}}}
\newcommand{\ad}{{\mathop{\rm ad}}}
\newcommand{\Ad}{{\mathop{\rm Ad}}}
\newcommand{\alg}{{\mathop{\rm alg}}}
\newcommand{\clalg}{{\mathop{\overline{\rm alg}}}}
\newcommand{\Aut}{\mathop{\rm Aut}}
\newcommand{\arccot}{\mathop{\rm arccot}}
\newcommand{\capp}{{\mathop{\rm cap}}}
\newcommand{\rcapp}{{\mathop{\rm rcap}}}
\newcommand{\diam}{\mathop{\rm diam}}
\newcommand{\divv}{\mathop{\rm div}}
\newcommand{\codim}{\mathop{\rm codim}}
\newcommand{\RRe}{\mathop{\rm Re}}
\newcommand{\IIm}{\mathop{\rm Im}}
\newcommand{\Tr}{{\mathop{\rm Tr \,}}}
\newcommand{\Vol}{{\mathop{\rm Vol}}}
\newcommand{\card}{{\mathop{\rm card}}}
\newcommand{\rank}{\mathop{\rm rank}}
\newcommand{\supp}{\mathop{\rm supp}}
\newcommand{\sgn}{\mathop{\rm sgn}}
\newcommand{\essinf}{\mathop{\rm ess\,inf}}
\newcommand{\esssup}{\mathop{\rm ess\,sup}}
\newcommand{\Int}{\mathop{\rm Int}}
\newcommand{\lcm}{\mathop{\rm lcm}}
\newcommand{\loc}{{\rm loc}}
\newcommand{\HS}{{\rm HS}}
\newcommand{\Trn}{{\rm Tr}}
\newcommand{\n}{{\rm N}}
\newcommand{\WOT}{{\rm WOT}}

\newcommand{\at}{@}

\newcommand{\mod}{\mathop{\rm mod}}
\newcommand{\spann}{\mathop{\rm span}}
\newcommand{\one}{\mathds{1}}

\hyphenation{groups}
\hyphenation{unitary}

\newcommand{\tfrac}[2]{{\textstyle \frac{#1}{#2}}}

\newcommand{\ca}{{\cal A}}
\newcommand{\cb}{{\cal B}}
\newcommand{\cc}{{\cal C}}
\newcommand{\cd}{{\cal D}}
\newcommand{\ce}{{\cal E}}
\newcommand{\cf}{{\cal F}}
\newcommand{\ch}{{\cal H}}
\newcommand{\chs}{{\cal HS}}
\newcommand{\ci}{{\cal I}}
\newcommand{\ck}{{\cal K}}
\newcommand{\cl}{{\cal L}}
\newcommand{\cm}{{\cal M}}
\newcommand{\cn}{{\cal N}}
\newcommand{\co}{{\cal O}}
\newcommand{\cs}{{\cal S}}
\newcommand{\ct}{{\cal T}}
\newcommand{\cx}{{\cal X}}
\newcommand{\cy}{{\cal Y}}
\newcommand{\cz}{{\cal Z}}

\newlength{\hightcharacter}
\newlength{\widthcharacter}
\newcommand{\covsup}[1]{\settowidth{\widthcharacter}{$#1$}\addtolength{\widthcharacter}{-0.15em}\settoheight{\hightcharacter}{$#1$}\addtolength{\hightcharacter}{0.1ex}#1\raisebox{\hightcharacter}[0pt][0pt]{\makebox[0pt]{\hspace{-\widthcharacter}$\scriptstyle\circ$}}}
\newcommand{\cov}[1]{\settowidth{\widthcharacter}{$#1$}\addtolength{\widthcharacter}{-0.15em}\settoheight{\hightcharacter}{$#1$}\addtolength{\hightcharacter}{0.1ex}#1\raisebox{\hightcharacter}{\makebox[0pt]{\hspace{-\widthcharacter}$\scriptstyle\circ$}}}
\newcommand{\scov}[1]{\settowidth{\widthcharacter}{$#1$}\addtolength{\widthcharacter}{-0.15em}\settoheight{\hightcharacter}{$#1$}\addtolength{\hightcharacter}{0.1ex}#1\raisebox{0.7\hightcharacter}{\makebox[0pt]{\hspace{-\widthcharacter}$\scriptstyle\circ$}}}

\thispagestyle{empty}

\vspace*{1cm}
\begin{center}
{\Large\bf The Dirichlet-to-Neumann operator \\[5pt]
via hidden compactness} \\[4mm]

\large W. Arendt$^1$, A.F.M. ter Elst$^2$, J.B. Kennedy$^1$ and M. Sauter$^2$

\end{center}

\vspace{4mm}

\begin{center}
{\bf Abstract}
\end{center}

\begin{list}{}{\leftmargin=1.8cm \rightmargin=1.8cm \listparindent=10mm 
   \parsep=0pt}
\item
We show that to each symmetric elliptic operator of the form
\begin{displaymath}
	\mathcal{A} = - \sum \partial_k \, a_{kl} \, \partial_l + c
\end{displaymath}
on a bounded Lipschitz domain $\Omega \subset \mathbb{R}^d$ one can 
associate a self-adjoint Dirichlet-to-Neumann operator on $L_2(\partial 
\Omega)$, which may be multi-valued if $0$ is in the Dirichlet spectrum of 
$\mathcal{A}$.
To overcome the lack of coerciveness in this case, we employ a new version 
of the Lax--Milgram lemma based on an indirect ellipticity property that 
we call hidden compactness.
We then establish uniform resolvent convergence of
 a sequence of Dirichlet-to-Neumann operators whenever 
their coefficients converge uniformly
and the second-order limit operator in $L_2(\Omega)$ has the 
unique continuation property.
We also consider semigroup convergence.

\end{list}

\vspace{6mm}
\noindent
April 2013.

\vspace{3mm}

\vspace{3mm}

\vspace{6mm}

\noindent
{\bf Home institutions:}    \\[3mm]
\begin{tabular}{@{}cl@{\hspace{10mm}}cl}
1. & Institute of Applied Analysis & 
2. & Department of Mathematics  \\
& University of Ulm & 
  & University of Auckland  \\
& Helmholtzstr.\ 18  &
  & Private bag 92019 \\
& D-89069 Ulm  &
   & Auckland 1142 \\ 
& Germany  &
  & New Zealand \\[8mm]
& wolfgang.arendt@uni-ulm.de & & terelst@math.auckland.ac.nz  \\
& james.kennedy@uni-ulm.de & & m.sauter@math.auckland.ac.nz  
\end{tabular}

\newpage
\setcounter{page}{1}

\section{Introduction} \label{Sdtonl1}

Let $\Delta^D$ be the Dirichlet Laplacian on a bounded Lipschitz domain
$\Omega \subset \Ri^d$ with boundary $\Gamma = \partial \Omega$.
Then for all $\lambda \in \Ri \setminus \sigma(-\Delta^D)$
one can define the Dirichlet-to-Neumann 
operator $D_\lambda$ as a self-adjoint operator on $L_2(\Gamma)$ 
which can be described by its graph
\begin{eqnarray*}
D_\lambda
= \{ (g,h) \in L_2(\Gamma) \times L_2(\Gamma) 
& : & 
   \mbox{there exists a } u \in H^1(\Omega) \mbox{ such that } \nonumber  \\*
& &   - \Delta u = \lambda \, u \mbox{ weakly in $\Omega$, } 
 u|_\Gamma = g
   \mbox{ and } \partial_\nu u = h \} 
.
\end{eqnarray*}
Here $u|_\Gamma$ is the trace of $u$ on $\Gamma$ and $\partial_\nu u$ 
its outer normal derivative, and we identify the operator $D_\lambda$ 
 with its graph  in a natural way.
(See, e.g., \cite{AE3, AE2, ArM2, BehR, Daners4, GesM}, 
and the references therein.)

What is perhaps less well known is that it is still possible to give meaning to 
$D_\lambda$ if $\lambda \in \sigma(-\Delta^D)$.
In this case there is now a nontrivial solution to $-\Delta u = \lambda \, u$ 
in $\Omega$ with $u|_\Gamma = 0$.
For simplicity assume that $\Omega$ has a $C^2$-boundary, 
so that this solution $u \in W^{2,2}(\Omega)$.
Then for each $g \in D(D_\lambda)$, the 
domain of $D_\lambda$, there is no longer a unique $h \in L_2(\Gamma)$ 
for which $D_\lambda g = h$, since $h + \partial_\nu u$ is 
obviously also a solution.

However, if we consider $D_\lambda$ as a graph (which we will do throughout 
the paper), then $D_\lambda$ becomes the graph of a possibly 
multi-valued operator if $\lambda \in \sigma(-\Delta^D)$. 
In order to avoid confusion, we will henceforth always use the term `graph' to mean 
`multi-valued operator', reserving `operator' for the single-valued type.
It was shown in \cite{ArM2} that the graph $D_\lambda$ is in fact self-adjoint 
(see Section~\ref{Sdtonl3.1} for the precise definitions), which is a consequence 
of the range condition $R(D_\lambda + isI) = L_2(\Gamma)$ being satisfied
for all $s \in \Ri \setminus \{ 0 \} $.
In order to establish this, one cannot use the usual form methods, since 
coerciveness (and even ellipticity) of the associated form are lost.
In \cite{ArM2} an alternative argument based on a Galerkin approximation method
given by Gr\'egoire, N\'ed\'elec and Planchard \cite{GNP} was used.

The purpose of this paper is to develop a form method which 
can be used in the above setting, whose point of departure may be found in the 
framework introduced in \cite{AE3, AE2}.
This will allow us not only to give an alternative proof that $D_\lambda$ is a 
self-adjoint graph, but also to establish various other properties of more 
general `Dirichlet-to-Neumann graphs'.

The key component of our work is a new argument which we call `hidden compactness'.
It establishes the Fredholm alternative (injectivity implies invertibility) for an 
operator defined by a sesquilinear and continuous but 
non-coercive form $\gota \colon V \times V \to \Ci$ under 
the assumption that $\gota$ is `compactly elliptic', that is, that there exists 
another Hilbert space $\widetilde H$ and a compact map $\tilde j \in \cl(V,\widetilde H)$ 
such that $\gota$ is $\tilde j$-elliptic (see Lemma~\ref{ldtonl223}).
The space $\widetilde H$ and map $\tilde j$ may be essentially arbitrary, 
provided only that the compact ellipticity criterion is satisfied, and do not enter into 
the theory in any other way; hence the `hiddenness' (and the tildes).
This result, which we regard as a `Fredholm--Lax--Milgram lemma', contains the 
classical Fredholm alternative as a special case and may be used as a substitute 
for the usual Lax--Milgram lemma, allowing us to develop a general theory of 
Dirichlet-to-Neumann graphs.
Although we will only be considering graphs, we wish to emphasize that 
this Fredholm--Lax--Milgram lemma and its application
are new, and possibly of general interest, even 
in the case of (single-valued) operators, as an addition to the general 
corpus of available form-theoretic tools.

In Section~\ref{Sdtonl2.1} we introduce the motivating example of a
Dirichlet-to-Neumann graph which will be of especial 
interest to us, and to which we will repeatedly return throughout the paper.
In Section~\ref{Sdtonl3.1} we introduce a number of essential definitions and 
basic results in the study of self-adjoint graphs in order to fix notation and 
terminology, and to keep the paper more self-contained.
In Section~\ref{Sdtonl2222} we formally introduce the notion of hidden compactness, 
give our Fredholm--Lax--Milgram lemma, and use it to prove, among other properties, 
that compact ellipticity of the symmetric form $\gota$ implies 
self-adjointness of the associated Dirichlet-to-Neumann graph $A$ 
(Theorem~\ref{tdtonl202}), as well as the surprising result that $A$ is always 
bounded below (Theorem~\ref{tdtonl210.5}).
We also characterize the single-valued part of $A$ (Proposition~\ref{pdtonl210.7}) 
and cast our results in the setting of the (concrete) Dirichlet-to-Neumann 
graph from Section~\ref{Sdtonl2.1}.

The other main topic of interest of the paper, which is the focus of 
Sections~\ref{Sdtonl5}--\ref{Sdtonl7}, is the study of `approximation of graphs', 
that is, under what conditions one can expect convergence of the resolvents and 
semigroups associated with a sequence of Dirichlet-to-Neumann graphs 
$(A_n)_{n \in \Ni}$.
If $A$ is a self-adjoint graph, then for all $s \in \Ri \setminus\{0\}$ the resolvent 
$(A + i \, s \, I)^{-1}$ is a single-valued bounded operator on $L_2(\Gamma)$.
In Section~\ref{Sdtonl5} we give a useful and natural criterion on a sequence 
$(\gota_n)_{n \in \Ni}$ of forms converging weakly to another form $\gota$, 
which implies that $\lim (A_n + i \, s \, I)^{-1} = (A + i \, s \, I)^{-1}$
strongly, where $A_n$ and $A$ are the associated Dirichlet-to-Neumann graphs 
(see Theorems~\ref{tdtonl502.3} and \ref{tdtonl502}, the latter being arguably 
the deepest abstract result in the paper).
In Theorems~\ref{tdtonl509.5} 
and \ref{tdtonl509} we give an analogous criterion on the forms under which the 
graphs $A_n$ are uniformly bounded below.
In fact, it turns out that this property is independent of 
convergence in the strong (even uniform) resolvent sense; 
it seems that two different subspaces of $V$ associated with $\gota$, 
which we denote by $W(\gota)$ and $V(\gota)$ 
(defined in Section~\ref{Sdtonl2222}), emerge naturally when determining 
resolvent convergence and uniform lower boundedness of the $A_n$, 
respectively.
We also consider strong convergence of the associated 
semigroups in Section~\ref{tdtonl6}.

As an application to our specific operator/graph $D_\lambda$ we obtain 
$\lim_{\lambda \to \lambda_0} (D_\lambda + i \, s \, I)^{-1} 
      = (D_{\lambda_0} + i \, s \, I)^{-1}$
for all $\lambda_0 \in \Ri$, regardless of whether or not $\lambda_0\in\sigma(- \Delta^D)$.
This is the subject of Section~\ref{Sdtonl7} (see Theorem~\ref{tdtonl703}), where 
we also prove a similar statement for the corresponding semigroups 
(Theorem~\ref{tdtonl705}).

Our criterion is also applicable in the far more general setting of a sequence of 
second-order elliptic operators with real symmetric bounded measurable coefficients.
If the relevant coefficients converge uniformly, then the results from 
Section~\ref{Sdtonl5} imply that the associated Dirichlet-to-Neumann graphs converge 
uniformly in the resolvent sense, if the limit operator in $L_2(\Omega)$ satisfies the 
interesting additional hypothesis that it possesses the unique continuation property 
(see Theorem~\ref{tdtonl707}).
This latter property has received much attention in the literature.
It is known to hold, for example, in two dimensions (Schulz \cite{Schulz}), 
or in higher dimensions if the coefficients are Lipschitz continuous 
(cf.~Kurata \cite{Kur}), but not in general if the coefficients are only H\"older 
continuous (see \cite{Fil1}).

Finally, in Section~\ref{Sdtonl4}, we consider $m$-accretive graphs.
We prove that 
if the form $\gota$ is accretive and compactly elliptic, 
then the associated Dirichlet-to-Neumann graph is $m$-accretive 
(Theorem~\ref{tdtonl430}).

\section{The basic example} \label{Sdtonl2.1}

Throughout this paper we consider a basic example, namely the 
Dirichlet-to-Neumann graph $D_m$ associated with $- \Delta + m$, 
which arises naturally in the context of what may be thought of as the `classical' 
Dirichlet-to-Neumann operator.
Here we shall explain this basic example in more detail and only afterwards introduce 
the abstract tools which allow us to treat this and also much more general examples.
In several instances we shall explain the abstract notions and results in 
terms of this concrete example.

Let $\Omega \subset \Ri^d$ be a Lipschitz domain, i.e.\ an open bounded
non-empty set such that for all $z \in \partial \Omega$ there exists an $r > 0$
such that $B(z,r) \cap \partial \Omega$ is a Lipschitz graph with 
$B(z, r) \cap \Omega$ on one side.
On the boundary $\Gamma = \partial \Omega$ we consider the $(d-1)$-dimensional 
Hausdorff measure (the surface measure).
The space $L_2(\Gamma)$ is formed with respect to this measure.
We denote by 
\[
H^1(\Omega)
= \{ u \in L_2(\Omega) : \partial_j u \in L_2(\Omega) \mbox{ for all } 
      j \in \{ 1,\ldots,d \} \}
\]
the first Sobolev space with norm
$\|u\|_{H^1(\Omega)}^2 = \|u\|_{L_2(\Omega)}^2 + \sum_{j = 1}^d \|\partial_j u\|_{L_2(\Omega)}^2$.
The space $H^1_0(\Omega)$ is the closure of the space of test functions 
$\cd(\Omega) = C_c^\infty (\Omega)$ in $H^1(\Omega)$.
There exists a unique bounded operator
$\Tr \colon H^1(\Omega) \to L_2(\Gamma)$, called the {\bf trace operator}, 
such that 
\[
\Tr u = u|_\Gamma
\]
for all $u \in H^1(\Omega) \cap C(\overline \Omega)$.

Next we define the normal derivative in a weak form by requiring 
Green's formula to be valid.

\begin{definition} \label{ddtonl240}
Let $u \in H^1(\Omega)$ be such that $\Delta u \in L_2(\Omega)$ in the 
sense of distributions.
If $h \in L_2(\Gamma)$ and 
\[
\int_\Omega (\Delta u) \, \overline v
   + \int_\Omega \nabla u \cdot \overline{\nabla v}
= \int_\Gamma h \, \overline{\Tr v}
\]
for all $v \in H^1(\Omega)$, then we call $h$ the {\bf normal derivative}
of $u$ and write $\partial_\nu u = h$.
(The element $h$ is obviously unique if it exists.)
We write $\partial_\nu u \in L_2(\Gamma)$ if there exists an $h \in L_2(\Gamma)$
such that $\partial_\nu u = h$.
\end{definition}

Let $m \in L_\infty(\Omega,\Ri)$.
We denote by $- \Delta^D + m$ the realization of $- \Delta + m$ with 
Dirichlet boundary conditions, i.e.\
\[
D(- \Delta^D + m)
= \{ u \in H^1_0(\Omega) : \Delta u \in L_2(\Omega) \}
\]
and 
\[
(- \Delta^D + m) u = - \Delta u + m \, u
\]
for all $u \in D(- \Delta^D + m)$.
This is a self-adjoint operator with compact resolvent.

Now we introduce our basic example.

\begin{definition} \label{ddtonl241}
Given $m \in L_\infty(\Omega,\Ri)$, the 
{\bf Dirichlet-to-Neumann graph} $D_m$
is defined by
\begin{eqnarray*}
D_m
= \{ (g,h) \in L_2(\Gamma) \times L_2(\Gamma) 
& : & 
   \mbox{there exists a } u \in H^1(\Omega) \mbox{ such that }  \\*
& &   (- \Delta + m) u = 0, \;
 \Tr u = g
   \mbox{ and } \partial_\nu u = h \} 
.
\end{eqnarray*}
\end{definition}

If $0 \not\in \sigma(- \Delta^D + m)$, then $D_m$ is 
the graph of a (single-valued) self-adjoint operator,
but in general $D_m$ might be multi-valued, i.e.\
there exists an $h \in L_2(\Gamma)$ with $(0,h) \in D_m$ but $h \neq 0$.
Nevertheless, $D_m$ is a self-adjoint graph.
We will explain this notion in the following section.

\section{Self-adjoint graphs} \label{Sdtonl3.1}

Multi-valued operators play an important role in non-linear analysis
(see Br\'ezis \cite{Bre2} and Showalter \cite{Sho}).
As mentioned in the introduction, 
we will use the term {\em graph} instead of `(possibly) multi-valued operator' 
and always consider linear graphs in this paper.
We reserve the term {\em operator} for single-valued operators.
In this section we give an expository account of some basic 
properties of graphs.

Let $H$ be a complex Hilbert space.
A {\bf graph} is a subspace of $H \times H$.
Let $A$ be a graph.
Then for each $x \in H$ we define $A(x)$ as the {\em set}
\[
A(x) = \{ y \in H : (x,y) \in A \}.
\]
Moreover, we set 
\begin{eqnarray*}
D(A) & = & \{ x \in H : \mbox{there exists an } y \in H \mbox{ such that } (x,y) \in A \} \mbox{ and }  \\
R(A) & = & \{ y \in H : \mbox{there exists a } x \in H \mbox{ such that } (x,y) \in A \} ,
\end{eqnarray*}
to be the {\bf domain} and {\bf range} of the graph~$A$, respectively.
The graph $A$ is called {\bf surjective} if $R(A) = H$.
We denote the reflection of $A$ in the diagonal of $H \times H$ by $A^\dagger$.
So 
\[
A^\dagger = \{ (y,x) : (x,y) \in A \} 
.  \]
We call $A$  
{\bf single-valued} if $A(x)$ has at most one element for all $x \in H$.
This is equivalent to $A(0) = \{ 0 \} $.
We call $A$ {\bf invertible} if $A^\dagger$ is single-valued, $A$ is surjective and 
$A$ is closed.
If the {\em graph} $A$ is invertible then one can define the 
{\em operator} $A^{-1} \colon H \to H$
by $A^{-1} y = x$ if $(x,y) \in A$.
It follows from the closed graph theorem that $A^{-1}$ is a bounded operator.
If $\lambda \in \Ci$ then define the graph $A + \lambda \, I$ by
\[
A + \lambda \, I
= \{ (x, y + \lambda \, x) : (x,y) \in A \}
.  \]
Define the {\bf resolvent set} $\rho(A)$ by
\[
\rho(A) = \{ \lambda \in \Ci : \mbox{the graph } A - \lambda \, I \mbox{ is invertible} \}
.  \]
It is easy to verify the resolvent identity 
\begin{equation}
(A - \lambda \, I)^{-1} - (A - \mu \, I)^{-1}
= (\lambda - \mu) \, (A - \lambda \, I)^{-1} \, (A - \mu \, I)^{-1}
\label{eSdtonl3.1;1}
\end{equation}
for all $\lambda,\mu \in \rho(A)$.
We say that the graph $A$ has {\bf compact resolvent} if there exists a 
$\lambda \in \rho(A)$ such that 
$(A - \lambda \, I)^{-1}$ is compact.
By (\ref{eSdtonl3.1;1}) this is equivalent to 
$(A - \lambda \, I)^{-1}$ being compact for all $\lambda \in \rho(A)$.

We call the graph $A$ {\bf symmetric} if $(x,y)_H \in \Ri$ for all $(x,y) \in A$.
The graph $A$ is called {\bf self-adjoint} if $A$ is symmetric and 
for all $s \in \Ri \setminus \{ 0 \} $ the graph $A + i \, s \, I$ is surjective.
Finally, a self-adjoint graph is called {\bf bounded below} if
there exists an $\omega \in \Ri$ such that 
\[
(x,y)_H + \omega \, \|x\|_H^2 \geq 0
\]
for all $(x,y) \in A$.
If $\omega$ can be taken as $0$, then $A$ is called {\bf positive}.

\begin{lemma} \label{ldtonl340}
Let $A$ be a self-adjoint graph in a Hilbert space $H$.
Then $i \, \Ri \setminus \{ 0 \} \subset \rho(A)$.
Moreover, $\|(A + i \, s \, I)^{-1}\| \leq \frac{1}{|s|}$
for all $s \in \Ri \setminus \{ 0 \} $.
\end{lemma}
\proof\
Let $s \in \Ri \setminus \{ 0 \} $.
Let $x \in (A + i \, s \, I)^\dagger(0)$.
Then $(x,0) \in A + i \, s \, I$.
Hence $(x, -i \, s \, x) \in A$.
Since $A$ is real one deduces that 
$i \, s \, \|x\|_H^2 = (x, -i \, s \, x)_H \in \Ri$.
So $x = 0$.

If $(x,y) \in A + i \, s \, I$, then 
$(x, y - i \, s \, x) \in A$.
Therefore $(x,y)_H + i \, s \, \|x\|_H^2 = (x, y - i \, s \, x)_H \in \Ri$ and 
$|s| \, \|x\|_H^2 = |\IIm (x,y)_H| \leq \|x\|_H \, \|y\|_H$.
Consequently, $|s| \, \|x\|_H \leq \|y\|_H$.
This implies that $A + i \, s \, I$ is closed and hence invertible.
Then the norm estimate is obvious.\hfill$\Box$

\vertspace

We list some properties of self-adjoint graphs.

\begin{prop} \label{pdtonl341}
Let $A$ be a self-adjoint graph.
\begin{tabel}
\item  \label{pdtonl341-1}
The set $A$ is closed in $H \times H$.
\item  \label{pdtonl341-2}
If $s \in \Ri \setminus \{ 0 \} $, then 
$( (A + i \, s \, I)^{-1} )^* = (A - i \, s \, I)^{-1}$.
\item  \label{pdtonl341-3}
If $(x,y), (u,v) \in A$, then 
$(u,y)_H = (v,x)_H$.
\item  \label{pdtonl341-4}
$A(0) = \ker (A + i \, s \, I)^{-1}$
for all $s \in \Ri \setminus \{ 0 \} $.
\item  \label{pdtonl341-5}
The graph $A$ is single valued if and only if $(A + i \, s \, I)^{-1}$
is injective for all $($or for one$)$ $s \in \Ri \setminus \{ 0 \} $.
\end{tabel}
\end{prop}

Note that \ref{pdtonl341-3} is the condition
$(u,B x)_H = (Bu, x)_H$ for all $x,u \in D(B)$ if $A$ is the graph of a 
self-adjoint operator $B$ in $H$.

Let $A$ be a self-adjoint graph. 
If $(x,y) \in A$ then $(x,y + y') \in A$ for all $y' \in A(0)$.
Therefore $(x,y + y')_H \in \Ri$ for all $y' \in A(0)$, which implies that 
$x \in A(0)^\perp$. 
So $D(A) \subset A(0)^\perp$.
Define the  operator $A^\circ$ in $A(0)^\perp$ by 
$D(A^\circ) = D(A)$ and 
\[
A^\circ x = y
\]
where $y \in A(0)^\perp$ is the unique element such that $(x,y) \in A$.
We call $A^\circ$ the {\bf single-valued part} of $A$.

\begin{prop} \label{pdton342}
Let $A$ be a self-adjoint graph.
\begin{tabel}
\item \label{pdton342-1}
The single-valued part $A^\circ$ of $A$ is a self-adjoint operator.
In particular, $D(A^\circ)$ is dense in $A(0)^\perp$.
\item \label{pdton342-2}
The graph $A$ is bounded below if and only if the single-valued part
$A^\circ$ is bounded below. 
\item \label{pdton342-3}
If $s \in \Ri \setminus \{ 0 \} $, then 
$(A + i \, s \, I)^{-1} = (A^\circ + i \, s \, I)^{-1} \oplus 0$,
where the decomposition is with respect to the decomposition
$H = A(0)^\perp \oplus A(0)$.
\end{tabel}
\end{prop}
\proof\
Clearly $(x,A^\circ x) \in A$, so $(x,A^\circ x)_H \in \Ri$ for all $x \in D(A^\circ)$.
Therefore $A^\circ$ is symmetric.

Next let $s \in \Ri \setminus \{ 0 \} $.
Let $y \in A(0)^\perp$.
Since $A$ is a self-adjoint graph, there exists a $x \in H$ such that 
$(x,y) \in A + i \, s \, I$.
Then $x \in D(A) = D(A^\circ)$ 
and $(A^\circ + i \, s \, I) x = y$.
Hence $A^\circ + i \, s \, I$ is 
surjective. 
Therefore $A^\circ$ is self-adjoint.
In particular, $A^\circ$ is densely defined 
and $D(A^\circ)$ is dense in $A(0)^\perp$.

The other statements are easy.\hfill$\Box$

\vertspace

The following converse of Proposition~\ref{pdton342} is easy to see. 

\begin{prop} \label{pdton343}
Let $H_1$ be a closed subspace of $H$ and let $B$ be a self-adjoint 
operator in $H_1$.
Define 
\[
A = \{ (x,y + Bx) : x \in D(B) \mbox{ and } y \in H_1^\perp \}
.  \]
Then $A$ is a self-adjoint graph and $A^\circ = B$.
\end{prop}

\section{Self-adjointness via hidden compactness} \label{Sdtonl2222}

The aim of this section is to give a criterion, which we call hidden 
compactness, to prove that a graph is self-adjoint.
First we introduce some notation and terminology.

Let $V$ be a complex Hilbert space and let 
$\gota \colon V \times V \to \Ci$ be a continuous sesquilinear form.
The form $\gota$ is called {\bf coercive} if there exists a $\mu > 0$ such that 
\[
\RRe \gota(u) \geq \mu \, \|u\|_V^2
\]
for all $u \in V$, where $\gota(u) = \gota(u,u)$.
Given a Hilbert space $H$ and $j \in \cl(V,H)$, we recall
from \cite{AE2} that the form 
$\gota$ is called {\bf $j$-elliptic} if there are 
$\omega \in \Ri$ and $\mu > 0$ such that 
\[
\RRe \gota(u) + \omega \, \|j(u)\|_H^2
\geq \mu \, \|u\|_V^2
\]
for all $u \in V$.
If $j$ is the inclusion of $V$ into $H$, then we also say that $\gota$ is 
{\bf $H$-elliptic} if $\gota$ is $j$-elliptic.
Next we introduce the following expression, which is new.
We say the form $\gota$ is {\bf compactly elliptic}
if there exists a Hilbert space $\widetilde H$ and a compact 
$\tilde j \in \cl(V,\widetilde H)$ such that $\gota$ is $\tilde j$-elliptic.
Clearly each coercive form is compactly elliptic.
In the next lemma, which we call the {\bf Fredholm--Lax--Milgram lemma},
the coerciveness condition in the original Lax--Milgram lemma
is replaced by compact ellipticity and an injectivity hypothesis.
Thus the hypothesis is a kind of hidden compactness, which will be central to 
establish self-adjointness and lower boundedness of our 
Dirichlet-to-Neumann graphs.
However, the map $\tilde j$ and the space $\widetilde H$ surprisingly serve no 
further purpose in the development of the subsequent general theory, and are 
therefore marked with tildes throughout to prevent confusion with other maps 
and spaces.

\begin{lemma} \label{ldtonl223}
Let $V$ be a Hilbert space and $\gota \colon V \times V \to \Ci$ 
a compactly elliptic continuous sesquilinear form.
Define the operator $\ca \colon V \to V'$ by 
\[
\gota(u,v) = (\ca u,v)_{V' \times V}
.  \]
Suppose that $\ca$ is injective.
Then $\ca$ is invertible.
\end{lemma}
\proof\
By assumption there exist a Hilbert space $\widetilde H$, a compact 
$\tilde j \in \cl(V,\widetilde H)$ and $\mu > 0$ such that 
$\RRe \gota(u) + \|\tilde j(u)\|_{\widetilde H}^2 \geq \mu \, \|u\|_V^2$ for all $u \in V$.
There exists a unique $T_0 \in \cl(V)$ such that $\gota(u,v) = (T_0 u,v)_V$ 
for all $u,v \in V$.
Define $T = T_0 + K$, where $K = {\tilde j}^* \, \tilde j$ is compact by assumption.
Then $\RRe (T u, u)_V \geq \mu \, \|u\|_V^2$ and $\|T u\|_V \geq \mu \, \|u\|_V$
for all $u \in V$.
Hence $T$ is injective and $T$ has closed range.
Similarly $T^*$ is injective.
Therefore $T$ is invertible. 
Since $T_0 = T(I - T^{-1} \, K)$ is injective and $T^{-1} \, K$ is 
compact, the operator $T_0$ is invertible by the Fredholm alternative
for $(I - T^{-1} \, K)$.
This is equivalent to $\ca$ being invertible.\hfill$\Box$

\begin{remarkn} \label{cdtonl223.5}
Lemma~\ref{ldtonl223} contains the classical Fredholm alternative as a 
special case.
In fact, let $\widetilde H$ be a Hilbert space and let $K \in \cl(\widetilde H)$ 
be compact.
Suppose that $I + K$ is injective.
Then $I + K$ is surjective.
Just choose $V = \widetilde H$,  $\tilde j = K$ and  $\gota(u,v) = 
((I + K)u,v)_{\widetilde H}$ in Lemma~\ref{ldtonl223}.
\end{remarkn}

Not every continuous sesquilinear form is compactly elliptic, 
as the following simple example shows.

\begin{exam} \label{xdtonl223.6}
Let $V$ be a Hilbert space and $\gota \colon V \times V \to \Ci$ given 
by $\gota(u,v) = - (u,v)_V$.
Then $\gota$ is compactly elliptic if and only if $V$ is finite dimensional.
Indeed, if $\widetilde H$ is a Hilbert space, $\tilde j \in \cl(V,\widetilde H)$ is compact 
and $\gota$ is $\tilde j$-elliptic, then there exists an $\alpha > 0$ such that 
$\|\tilde j(u)\|_{\widetilde H}^2 \geq \alpha \, \|u\|_V^2$ for all $u \in V$.
Then there exists an $S \in \cl(\widetilde H,V)$ such that 
$S \circ \tilde j = I_V$.
So the identity operator $I_V$ on $V$ is compact.
\end{exam}

Compact ellipticity has some useful permanence properties.

\begin{prop} \label{pdtonl223.7}
Let $V$ be a Hilbert space and $\gota \colon V \times V \to \Ci$ be a 
continuous compactly elliptic form.
\begin{tabel}
\item \label{pdtonl223.7-1}
If $\gotb \colon V \times V \to \Ci$ is a compactly elliptic form, then 
so is $\gota + \gotb$.
\item \label{pdtonl223.7-2}
Let $K \in \cl(V)$ be compact.
Define $\gotb \colon V \times V \to \Ci$ by $\gotb(u,v) = \gota(u,v) + (Ku,v)_V$.
Then $\gotb$ is a compactly elliptic form.
\item \label{pdtonl223.7-3}
Let $V_1$ be a closed subspace of $V$. 
Then $\gota|_{V_1 \times V_1}$ is compactly elliptic.
\end{tabel}
\end{prop}
\proof\
By assumption there exist a Hilbert space $\widetilde H$, a compact 
$\tilde j \in \cl(V,\widetilde H)$ and $\mu > 0$ such that 
$\RRe \gota(u) + \|\tilde j(u)\|_{\widetilde H}^2 \geq \mu \, \|u\|_V^2$ for all $u \in V$.

`\ref{pdtonl223.7-1}'.
Since $\gotb$ is compactly elliptic, 
there exist a Hilbert space ${\widetilde H}_1$, a compact ${\tilde j}_1 \in 
\cl(V,{\widetilde H}_1)$ and $\mu_1 > 0$ such that 
$\RRe \gotb(u) + \|{\tilde j}_1(u)\|_{{\widetilde H}_1}^2 \geq \mu_1 \, 
\|u\|_V^2$ for all $u \in V$.
Choose ${\widetilde H}_2 = \widetilde H \oplus {\widetilde H}_1$ and ${\tilde j}_2 = 
\tilde j \oplus {\tilde j}_1$.
Then ${\tilde j}_2$ is compact and $\gota + \gotb$ is ${\tilde j}_2$-elliptic.

`\ref{pdtonl223.7-2}'.
Choose ${\widetilde H}_3 = \widetilde H \oplus V$ and define ${\tilde j}_3 \in 
\cl(V,{\widetilde H}_3)$ by ${\tilde j}_3(u) = (\tilde j(u), Ku)$.
Then ${\tilde j}_3$ is compact.
Let $u \in V$.
Then 
\[
|(Ku,u)_V| 
\leq \|Ku\|_V \, \|u\|_V
\leq \frac{\mu}{2} \, \|u\|_V^2 + \frac{1}{2 \mu} \|Ku\|_V^2
.  \]
Therefore 
\[
\RRe \gotb(u) + (1 + \frac{1}{2 \mu}) \, \|{\tilde j}_3(u)\|_{{\widetilde H}_3}^2
\geq \frac{\mu}{2} \, \|u\|_V^2
\]
and $\gotb$ is ${\tilde j}_3$-elliptic.

The last statement is easy.\hfill$\Box$

\vertspace

Given Hilbert spaces $V$ and $H$, a continuous form 
$\gota \colon V \times V \to \Ci$ and an operator $j \in \cl(V,H)$ 
we define the {\bf graph associated with $(\gota,j)$} in $H \times H$ by 
\begin{eqnarray*}
A = \{ (x,y) \in H \times H
& : & \mbox{there exists a } u \in V \mbox{ such that}  \\*
& & j(u) = x \mbox{ and } \gota(u,v) = (y,j(v))_H \mbox{ for all } v \in V \}
.
\end{eqnarray*}
We consider $A$ as an abstract Dirichlet-to-Neumann graph.
If $\gota$ is $j$-elliptic and $j(V)$ is dense in $H$, then $A$ is
the graph of a (single-valued) sectorial operator (see \cite{AE2} Theorem~2.1).
The following is the main result of this section.
It is a generation theorem where we replace $j$-ellipticity
by the condition that $\gota$ is compactly elliptic, i.e.\
we assume the existence of a compact operator $\tilde j \in \cl(V,\widetilde H)$
for which the form $\gota$ is $\tilde j$-elliptic.
We emphasize that the maps $j$ and $\tilde j$ are different in general.

We need the following subspace of~$V$.
For any form $\gota$ on $V$ and fixed $j \in \cl(V,H)$ define 
\[
W(\gota) = \{ u \in \ker j : \gota(u,v) = 0 \mbox{ for all } v \in V \}
.  \]
This space is always taken with respect to the map~$j$.
The map $\tilde j$, used in compact ellipticity, plays no role in the 
definition of $W(\gota)$ and the graph~$A$.

This space will play a decisive role later in Section~\ref{Sdtonl5}, 
but it will also be used in the proof of Theorem~\ref{tdtonl202}.
Note that if $(x,y) \in A$ and $u_0 \in V$ is such that $j(u_0) = x$ and 
$\gota(u_0,v) = (y,j(v))_H$ for all $v \in V$, then 
\begin{equation}
\{ u \in V : j(u) = x \mbox{ and } \gota(u,v) = (y,j(v))_H \mbox{ for all } v \in V \}
= u_0 + W(\gota)
.  
\label{etdton202;1}
\end{equation}
Therefore we call $W(\gota)$ the space of non-uniqueness.

\begin{thm} \label{tdtonl202}
Let $V$ and $H$ be Hilbert spaces.
Let $\gota \colon V \times V \to \Ci$ be a symmetric continuous sesquilinear form.
Further, let $j \in \cl(V,H)$.
Let $A$ be the graph associated with $(\gota,j)$.
If $\gota$ is compactly elliptic, then $A$ is a self-adjoint graph.
\end{thm}
\proof\
First suppose that $W(\gota) = \{ 0 \} $.

Let $(x,y) \in A$.
Let $u \in V$ be such that $j(u) = x$ and $\gota(u,v) = (y, j(v))_H$
for all $v \in V$.
Then 
\[
(x,y)_H 
= \overline{(y,x)_H}
= \overline{(y,j(u))_H}
= \overline{\gota(u)}
\in \Ri
.  \]
Next let $s \in \Ri \setminus \{ 0 \} $.
We shall show that $A + i \, s \, I$ is surjective.
Define the sesquilinear form $\gotb \colon V \times V \to \Ci$ by
\[
\gotb(u,v) = \gota(u,v) + i \, s \, (j(u), j(v))_H
.  \]
Since $j$ is continuous it follows that the form $\gotb$ is continuous.
Because $\gota$ is compactly elliptic, there exist a Hilbert space
$\widetilde H$, a compact $\tilde j \in \cl(V,\widetilde H)$ 
and $\mu  > 0$
such that $\gota(u) + \|\tilde j(u)\|_{\widetilde H}^2 \geq \mu \, \|u\|_V^2$
for all $u \in V$.
Then 
\[
\RRe \gotb(u) + \|\tilde j(u)\|_{\widetilde H}^2
= \gota(u) + \|\tilde j(u)\|_{\widetilde H}^2 
\geq \mu \, \|u\|_V^2
\]
for all $u \in V$.
Therefore $\gotb$ is $\tilde j$-elliptic.
Define $\cb \colon V \to V'$ by 
$(\cb u,v)_{V' \times V} = \gotb(u,v)$.
We show that $\cb$ is injective.

Let $u \in V$ and suppose that $\cb u = 0$.
Then 
\[
0 
= (\cb u, u)_{V' \times V} 
= \gotb(u) 
= \gota(u) + i \, s \, \|j(u)\|_H^2
.  \]
Since $\gota(u) \in \Ri$ and $s \in \Ri \setminus \{ 0 \} $ this 
implies that $j(u) = 0$.
Then for all $v \in V$ one has 
\[
0 
= (\cb u,v)_{V' \times V}
= \gotb(u,v)
= \gota(u,v) + i \, s \, (j(u), j(v))_H
= \gota(u,v)
.  \]
So $u \in W(\gota) = \{ 0 \} $ by assumption.
So $\cb$ is injective and therefore also surjective by the Fredholm--Lax--Milgram lemma,
Lemma~\ref{ldtonl223}.
Now let $y \in H$.
Define $\alpha \colon V \to \Ci$ by 
$\alpha(v) = (y, j(v))_H$.
Then $\alpha \in V'$ since $j$ is continuous.
Because $\cb$ is surjective, there exists a (unique) $u \in V$ such that
$\cb u = \alpha$.
Then for all $v \in V$ one has
\begin{eqnarray*}
(y, j(v))_H
& = & (\cb u,v)_{V' \times V}
= \gotb(u,v)
= \gota(u,v) + i \, s \, (j(u), j(v))_H
= \gota(u,v) + i \, s \, (x, j(v))_H
,
\end{eqnarray*}
where $x = j(u)$.
So $(x,y) \in A + i \, s \, I$.
This proves that $A$ is a self-adjoint graph if $W(\gota) = \{ 0 \} $.

Finally we drop the assumption that $W(\gota) = \{ 0 \} $.
Let $V_1 = W(\gota)^\perp$, where the orthogonal complement is in $V$.
Define $\gota_1 = \gota|_{V_1 \times V_1}$ and $j_1 = j|_{V_1}$.
Then $\gota_1$ is compactly elliptic by Proposition~\ref{pdtonl223.7}\ref{pdtonl223.7-3}.
Let $u \in W(\gota_1)$.
Then $u \in V_1$, $j(u) = 0$ and $\gota(u,v) = 0$ for all $v \in V_1$.
If $w \in W(\gota)$ then $\gota(w,u) = 0$ by definition of $W(\gota)$.
So $\gota(u,w) = \overline{\gota(w,u)} = 0$.
Hence by linearity $\gota(u,v) = 0$ for all $v \in V$.
Therefore $u \in W(\gota)$.
So $u \in W(\gota) \cap V_1 \subset W(\gota) \cap W(\gota)^\perp = \{ 0 \} $.
Thus $W(\gota_1) = \{ 0 \} $.
Let $A_1$ be the graph associated with $(\gota_1,j_1)$.
By the first part of the proof one deduces that $A_1$ is a self-adjoint graph.
In the next lemma we show that $A = A_1$.
Hence $A = A_1$ is a self-adjoint graph.\hfill$\Box$

\vertspace

To complete the proof of Theorem~\ref{tdtonl202} 
it remains to show the following general fact.

\begin{lemma} \label{ldtonl202.5}
Let $V$ and $H$ be Hilbert spaces and $\gota \colon V \times V \to \Ci$
a continuous symmetric sesquilinear form.
Let $j \in \cl(V,H)$.
Define $V_1 = W(\gota)^\perp$, $\gota_1 = \gota|_{V_1 \times V_1}$ and $j_1 = j|_{V_1}$.
Let $A$ and $A_1$ be the graphs associated with 
$(\gota,j)$ and $(\gota_1,j_1)$.
Then $A = A_1$.
\end{lemma}
\proof\
`$A_1 \subset A$'.
Let $x,y \in H$ and suppose that $(x,y) \in A_1$.
Then there exists a $u \in V_1$ such  that $j_1(u) = x$ and 
$\gota_1(u,v) = (y, j_1(v))_H$ for all $v \in V_1$.
Let $v \in V$.
Write $v = w + v_1$ with $w \in W(\gota)$ and $v_1 \in V_1$.
Then $j(w) = 0$.
Moreover, $\gota(u,w) = \overline{\gota(w,u)} = 0$.
So 
\[
\gota(u,v) 
= \gota(u,v_1)
= \gota_1(u,v_1)
= (y, j_1(v_1))_H
= (y, j(v_1))_H
= (y, j(v))_H
.  \]
Therefore $(x,y) \in A$.

`$A \subset A_1$'.
Let $x,y \in H$ and suppose that $(x,y) \in A$.
Let $u \in V$ be such that 
$j(u) = x$ and $\gota(u,v) = (y, j(v))_H$ for all $v \in V$.
Write $u = w + u_1$ with $w \in W(\gota)$ and $u_1 \in V_1$.
Then $\gota(w,v) = 0$ for all $v \in V$ and $j(w) = 0$.
So $j(u_1) = x$ and 
\[
\gota_1(u_1,v)
= \gota(u_1,v)
= \gota(u,v)
= (y,j(v))_H
\]
for all $v \in V_1$.
Therefore $(x,y) \in A_1$.\hfill$\Box$

\begin{remarkn} \label{rdtonl202.7}
Under the assumptions of Theorem~\ref{tdtonl202} the space $W(\gota)$ 
is finite dimensional.
Indeed, if $\widetilde H$ is a Hilbert space, $\tilde j \in \cl(V,\widetilde H)$ 
is compact and $\mu > 0$ are such that 
$\gota(u) + \|\tilde j(u)\|_{\widetilde H}^2 \geq \mu \, \|u\|_V^2$, then 
$\|\tilde j(u)\|_{\widetilde H}^2 \geq \mu \, \|u\|_V^2$ for all $u \in W(\gota)$.
Since $\tilde j|_{W(\gota)}$ is compact, the space $W(\gota)$ must be finite dimensional.
\end{remarkn}

In Theorem~\ref{tdtonl210.5} we shall prove that the self-adjoint graph $A$ is 
bounded below.
But first we prove that, as for  sectorial forms, the associated graph 
has compact resolvent if the map $j$ is compact.

\begin{prop} \label{pdtonl202.3}
Adopt the assumptions and notation of Theorem~{\rm \ref{tdtonl202}}.
In addition assume that the operator $j$ is compact.
Then $A$ has compact resolvent.
\end{prop}
\proof\
Using Lemma~\ref{ldtonl202.5} we may assume without loss of generality that 
$W(\gota) = \{ 0 \} $.
Let $s \in \Ri \setminus \{ 0 \} $.
Let $\cb$ be as in the proof of Theorem~\ref{tdtonl202}.
Then $\cb$ is invertible and $(A + i \, s \, I)^{-1} = j \circ \cb^{-1} \circ j^*$.
Since $j$ is compact, this resolvent is also compact.\hfill$\Box$

\vertspace

Before proving some additional properties in the situation 
of Theorem~\ref{tdtonl202}, we show that the basic example of Section~\ref{Sdtonl2.1}
is a self-adjoint graph with compact resolvent.

\begin{exam} \label{xdtonl450}
Let $\Omega \subset \Ri^d$ be a Lipschitz domain and let $m \in L_\infty(\Omega,\Ri)$.
Then 
\begin{eqnarray*}
D_m
= \{ (g,h) \in L_2(\Gamma) \times L_2(\Gamma) 
& : & 
   \mbox{there exists a } u \in H^1(\Omega) \mbox{ such that }  \\*
& &   (- \Delta + m) u = 0, \;
 \Tr u = g
   \mbox{ and } \partial_\nu u = h \} 
.
\end{eqnarray*}
is a self-adjoint graph with compact resolvent.

To see this, choose $H = L_2(\Gamma)$, $V = H^1(\Omega)$ and 
$j \colon H^1(\Omega) \to L_2(\Gamma)$ as the trace operator.
Define $\gota \colon V \times V \to \Ci$ by
\begin{equation}
\gota(u,v) 
= \int_\Omega \nabla u \cdot \overline{\nabla v}
   + \int_\Omega m \, u \, \overline v
.  \label{exdtonl450;1}
\end{equation}
Then $D_m$ is the graph associated with $(\gota,j)$.
This can be shown as follows.
Let $g,h \in L_2(\Gamma)$.
Suppose that $(g,h)$ is an element of the graph associated with $(\gota,j)$.
Then there exists a $u \in H^1(\Omega)$ such that 
$\Tr u = g$ and 
\begin{equation}
\int_\Omega \nabla u \cdot \overline{\nabla v}
   + \int_\Omega m \, u \, \overline v
= \int_\Gamma h \, \overline{\Tr v}
\label{exdtonl450;2}
\end{equation}
for all $v \in H^1(\Omega)$.
Taking $v \in \cd(\Omega)$ we see that $(-\Delta + m)u = 0$.
Replacing $m \, u$ by $\Delta u$ in (\ref{exdtonl450;2}) 
we deduce that $\partial_\nu u = h$.
Hence $(g,h) \in D_m$.
The converse inclusion is proved similarly.

We choose as $\tilde j$ the inclusion of $H^1(\Omega)$ into $L_2(\Omega)$.
This is a compact map and $\gota$ is $\tilde j$-elliptic.
Now it follows from Theorem~\ref{tdtonl202} that $D_m$ is a self-adjoint graph.
Since the trace operator is also compact, it follows from Proposition~\ref{pdtonl202.3}
that $D_m$ has compact resolvent.

It is not obvious that $D_m$ is lower bounded.
This follows from Theorem~\ref{tdtonl210.5} below, which needs further preparation.
\end{exam}

In order to prove that the graph of Theorem~\ref{tdtonl202} is 
bounded below, we need some reduction properties which are of 
independent interest.

Let $V$ and $H$ be Hilbert spaces, let
$\gota \colon V \times V \to \Ci$ be a continuous sesquilinear form and 
let $j \in \cl(V,H)$.
Further, let $V_1 \subset V$ be a closed subspace.
We define the restriction to the space $V_1$ by 
$\gota_1 = \gota|_{V_1 \times V_1}$ and $j_1 = j|_{V_1}$.
Let $A$ be the graph associated with $(\gota,j)$ and 
$A_1$ be the graph associated with $(\gota_1,j_1)$.
In general there is no relation between $A$ and $A_1$.
Even if one knows that $A \subset A_1$ or $A_1 \subset A$,
then it is still possible that the inclusion is strict.
But if $\gota$ is compactly elliptic, then so is $\gota_1$.
Hence both graphs are self-adjoint
and an inclusion implies equality.

\vertspace

We have to introduce one more space.
Let $V$ and $H$ be Hilbert spaces, let
$\gota \colon V \times V \to \Ci$ be a continuous sesquilinear form and 
let $j \in \cl(V,H)$.
Define 
\[
V(\gota) = \{ u \in V : \gota(u,v) = 0 \mbox{ for all } v \in \ker j \} 
.  \]
The space $V(\gota)$ will be used throughout this paper.
Its most important and immediate application is that we 
can consider the form $\gota$ restricted to this space 
and still obtain the same operator.
Thus we only need to consider the functions $u \in V$ `orthogonal' 
to $\ker j$ with respect to the form $\gota$.

\begin{prop} \label{pdtonl206}
Adopt the assumptions and notation of Theorem~{\rm \ref{tdtonl202}}.
Let $V_1 = V(\gota)$.
Define $\gota_1 = \gota|_{V_1 \times V_1}$ and $j_1 = j|_{V_1}$.
Let $A_1$ be the graph associated with $(\gota_1,j_1)$.
Then $A = A_1$.
Moreover, $V(\gota_1) = V(\gota)$.
\end{prop}
\proof\
Let $(x,y) \in A$.
Then there exists a $u \in V$ such that $j(u) = x$ and 
$\gota(u,v) = (y,j(v))_H$ for all $v \in V$.
In particular, $u \in V(\gota) = V_1$.
Moreover, for all $v \in V_1$ one has 
$\gota_1(u,v) = \gota(u,v) = (y,j(v))_H = (y,j_1(v))_H$.
Hence $(x,y) \in A_1$.
So $A \subset A_1$.
Since both graphs are self-adjoint one deduces that $A = A_1$.

Obviously $V(\gota_1) \subset V_1 = V(\gota)$.
Conversely, let $u \in V(\gota)$.
Then for all $v \in \ker j_1$ one has $v \in \ker j$
and hence $\gota_1(u,v) = \gota(u,v) = 0$.
So $u \in V(\gota_1)$ and $V(\gota) \subset V(\gota_1)$.
Therefore $V(\gota) = V(\gota_1)$.\hfill$\Box$

\begin{cor} \label{cdtonl207}
Adopt the assumptions and notation of Theorem~{\rm \ref{tdtonl202}}.
Let $V_1 = V(\gota) \cap (V(\gota) \cap \ker j)^\perp$.
Define $\gota_1 = \gota|_{V_1 \times V_1}$ and $j_1 = j|_{V_1}$.
Let $A_1$ be the graph associated with $(\gota_1,j_1)$.
Then $A = A_1$.
\end{cor}
\proof\
By Proposition~\ref{pdtonl206} we may without loss of generality 
assume that $V = V(\gota)$.
Let $(x,y) \in A_1$.
Then there exists a $u \in V_1$ such that 
$j_1(u) = x$ and $\gota_1(u,v) = (y,j_1(v))_H$ for all $v \in (\ker j)^\perp$.
Then $u \in V = V(\gota)$ and $\gota(u,v) = (y,j(v))_H$ for all $v \in (\ker j)^\perp$.
Also, if $v \in \ker j$ then $\gota(u,v) = 0 = (y,0)_H = (y,j(v))_H$.
So by linearity $\gota(u,v) = (y,j(v))_H$ for all $v \in V$.
Therefore $(x,y) \in A$ and $A_1 \subset A$.
By self-adjointness one deduces that $A = A_1$.\hspace*{5mm}\hfill$\Box$

\vertspace

We also need the following lemma, which shows that the hidden compactness 
argument does not cover a new situation if $j$ is injective.
Note that symmetry of the form $\gota$ is not required in the 
next lemma.

\begin{lemma} \label{ldtonl208}
Let $V$ be a Hilbert space and $\gota \colon V \times V \to \Ci$ be a 
compactly elliptic continuous form.
Further, let $H$ be a Hilbert space and $j \in \cl(V,H)$.
Suppose that $j$ is injective.
Then $\gota$ is $j$-elliptic.
\end{lemma}
\proof\
Because $\gota$ is compactly elliptic, there exist a Hilbert space
$\widetilde H$, a compact $\tilde j \in \cl(V,\widetilde H)$ 
and $\mu > 0$
such that 
\[
\RRe \gota(u) + \|\tilde j(u)\|_{\widetilde H}^2 
\geq \mu \, \|u\|_V^2
\]
for all $u \in V$.
Choose $\varepsilon = \frac{\mu}{2}$.
Since $\tilde j$ is compact and $j$ is injective, there exists a 
$c > 0$ such that 
$\|\tilde j(u)\|_{\widetilde H}^2 \leq \varepsilon \, \|u\|_V^2 + c \, \|j(u)\|_H^2$
for all $u \in V$.
Then 
\[
\RRe \gota(u) 
\geq (\mu - \varepsilon) \, \|u\|_V^2 - c \, \|j(u)\|_H^2
= \tfrac{1}{2} \, \mu \, \|u\|_V^2 - c \, \|j(u)\|_H^2
\]
for all $u \in V$ and $\gota$ is $j$-elliptic.\hfill$\Box$

\vertspace

Now we are able to prove lower boundedness.

\begin{thm} \label{tdtonl210.5}
Adopt the assumptions and notation of Theorem~{\rm \ref{tdtonl202}}.
Then $A$ is bounded below.
\end{thm}
\proof\
Let $V_1 = V(\gota) \cap (V(\gota) \cap \ker j)^\perp$.
Set $\gota_1 = \gota|_{V_1 \times V_1}$ and $j_1 = j|_{V_1}$.
Then $A$ is the graph associated with $(\gota_1,j_1)$
by Corollary~\ref{cdtonl207}.
But $j_1$ is injective.
So $\gota_1$ is $j_1$-elliptic by Lemma~\ref{ldtonl208}.
Hence there exists an $M \geq 0$ such that 
$\gota_1(u) + M \, \|j(u)\|_H^2 \geq 0$ for all $u \in V_1$.
Let $(x,y) \in A$.
Then there exists a $u \in V_1$ such that $j_1(u) = x$ and 
$\gota_1(u,v) = (y,j(v))_H$ for all $v \in V_1$.
Therefore 
\[
(y,x)_H
= (y,j(u))_H
= \gota_1(u)
\geq - M \|j(u)\|_H^2
= - M \, \|x\|_H^2
\]
and $A$ is bounded below.\hfill$\Box$

\vertspace

It is remarkable that the graph $A$ is bounded below, since the 
form in general is not bounded below (i.e.\ $j$-elliptic).
For example, consider a Lipschitz domain $\Omega$, 
let $\lambda_1^D$ and $\lambda_2^D$ denote the first and second 
eigenvalue of $- \Delta^D$.
Further, let 
$\lambda \in (\lambda_1^D,\lambda_2^D)$ and 
let $D_\lambda$ be the Dirichlet-to-Neumann graph 
associated with $(\gota,j)$, where the form $\gota \colon H^1(\Omega) \times H^1(\Omega) \to \Ci$ 
is given by 
\[
\gota(u,v) = \int_\Omega \nabla u \cdot \overline{\nabla v}
   - \lambda \int_\Omega u \, \overline v
\]
and $j = \Tr$ is the trace operator.
Then $D_\lambda$ is a self-adjoint graph which is bounded below 
by Theorem~\ref{tdtonl210.5}. 
But if $u$ is an eigenfunction of $- \Delta^D$ to the eigenvalue $\lambda_1^D$,
then 
\[
\gota(u) + M \, \|j(u)\|_{L_2(\Gamma)}^2
= \int_\Omega |\nabla u|^2 - \lambda \int_\Omega |u|^2
= (\lambda_1 - \lambda) \int_\Omega |u|^2
< 0
\]
for all $M \in \Ri$.
Thus the form associated with $D_\lambda$ is not bounded below.

\smallskip

Finally we determine the single-valued part of the 
self-adjoint graph in Theorem~{\rm \ref{tdtonl202}}.
For that we need one more lemma, which is also valid for non-symmetric 
forms.

\begin{lemma} \label{ldtonl209}
Let $V$ and $H$ be Hilbert spaces.
Let $\gota \colon V \times V \to \Ci$ be a continuous sesquilinear form
and $j \in \cl(V,H)$ be injective.
Suppose that $\gota$ is $j$-elliptic.
Let $A$ be the graph associated with $(\gota,j)$.
Then $\overline{j(V)} = A(0)^\perp$,
where the orthogonal complement is in $H$.
Moreover, the restriction 
$A|_{\overline{j(V)}}$ is a $($single-valued$)$ operator in $\overline{j(V)}$
which is $m$-sectorial.
\end{lemma}
\proof\
Let $y \in (j(V))^\perp$.
Then for all $v \in V$ one has 
$\gota(0,v) = 0 = (y,j(v))_H$.
Hence $(0,y) \in A$ and $y \in A(0)$.
Conversely, suppose that $y \in H$ and $(0,y) \in A$.
Then there exists a $u \in V$ such that $j(u) = 0$ and 
$\gota(u,v) = (y,j(v))_H$ for all $v \in V$.
Then $u = 0$ since $j$ is injective.
Moreover, $(y,j(v))_H = \gota(0,v) = 0$ for all $v \in V$.
So $y \in (j(V))^\perp$.
This proves the equality $\overline{j(V)} = A(0)^\perp$.

Let $j_1 \colon V \to H_1$ be the restriction of $j$, but with 
co-domain $H_1 = A(0)^\perp$.
Then $j_1 \in \cl(V,H_1)$ and $j_1(V)$ is dense in $H_1$.
By \cite{AE2} Theorem~2.1 one can associate a 
(single-valued) operator $A_1$ with $(\gota,j_1)$.
Then $A_1$ is $m$-sectorial.
It is straightforward to see that 
$G(A_1) = (H_1 \times H_1) \cap A$.
This proves the lemma.\hfill$\Box$

\vertspace

We are now able to characterize the single-valued part of the 
self-adjoint graph in Theorem~{\rm \ref{tdtonl202}}.

\begin{prop} \label{pdtonl210.7}
Adopt the assumptions and notation of Theorem~{\rm \ref{tdtonl202}}.
Let 
\[
H_1 = \overline{j(V(\gota))}
,  \]
where the closure is in $H$.
Then 
$A \cap (H_1 \times H_1)$ is the graph of the single-valued part of $A$.
\end{prop}
\proof\
Since $j(V(\gota)) = j(V(\gota) \cap (V(\gota) \cap \ker j)^\perp)$,
this follows from Corollary~\ref{cdtonl207} and Lemma~\ref{ldtonl209}.\hfill$\Box$

\section{Resolvent convergence} \label{Sdtonl5}

We now wish to investigate the convergence 
of a sequence of Dirichlet-to-Neumann graphs, where `convergence' is generally understood 
to be either of the associated resolvents or the semigroups.

In this section we consider the resolvent convergence 
$\lim_{n \to \infty} (A_n + i \, s \, I)^{-1}$ in various operator topologies.
Although we will generally consider only self-adjoint graphs $A_n$, 
our first result concerns resolvent convergence for arbitrary graphs.

\begin{prop} \label{pdtonl602}
Let $A,A_1,A_2,\ldots$ be graphs.
Let $\lambda,\mu \in \Ci$ and suppose that $\lambda,\mu \in \rho(A_n) \cap \rho(A)$
for all $n \in \Ni$.
Suppose that $\sup \|(A_n - \lambda \, I)^{-1}\| < \infty$
and $\sup \|(A_n - \mu \, I)^{-1}\| < \infty$.
Finally, suppose that 
$\lim_{n \to \infty} (A_n - \lambda \, I)^{-1} y = (A - \lambda \, I)^{-1} y$
for all $y \in H$.
Then $\lim_{n \to \infty} (A_n - \mu \, I)^{-1} y = (A - \mu \, I)^{-1} y$
for all $y \in H$.
\end{prop}
\proof\
This follows as in \cite{Kat1} Theorem~IV.2.25.\hfill$\Box$

\vertspace

Let $A,A_1,A_2,\ldots$ be graphs.
We say that $\lim_{n \to \infty} A_n = A$ in the {\bf strong resolvent sense}
if $\lim_{n \to \infty} (A_n - \lambda \, I)^{-1} = (A - \lambda \, I)^{-1}$
strongly for one (equivalently all) $\lambda \in \Ci$
with $\lambda \in \rho(A_n) \cap \rho(A)$ and  $\sup \|(A_n - \lambda \, I)^{-1}\| < \infty$.
If $A_n$ and $A$ are self-adjoint, this is equivalent to 
$\lim_{n \to \infty} (A_n + i \, s \, I)^{-1} = (A + i \, s \, I)^{-1}$
strongly for one (or all) $s \in \Ri \setminus \{ 0 \} $.

Throughout this section we fix Hilbert spaces $V$, $H$ and $\widetilde H$,
a continuous map $j \colon V \to H$ and a compact map $\tilde j \colon V \to \widetilde H$.
Further, for all $n \in \Ni$ let $\gota_n,\gota \colon V \times V \to \Ci$ be continuous
symmetric sesquilinear forms.
For all $n \in \Ni$ let $A_n$ be the graph associated with $(\gota_n,j)$ 
and let $A$ be the graph associated with $(\gota,j)$.

We say that the sequence $(\gota_n)_{n \in \Ni}$ is {\bf uniformly $\tilde j$-elliptic}
if there exist $\mu,\omega > 0$ such that
\begin{equation}
\gota_n(u) + \omega \, \|\tilde j(u)\|_{\widetilde H}^2
\geq \mu \, \|u\|_V^2
\label{eSdtonl5;1}
\end{equation}
for all $n \in \Ni$ and $u \in V$.
In addition, we say that $(\gota_n)_{n \in \Ni}$ {\bf converges weakly to} $\gota$ if 
\begin{equation}
\lim_{n \to \infty} \gota_n(u_n,v) = \gota(u,v)
\label{eSdtonl5;2}
\end{equation}
for all $v \in V$ and $u,u_1,u_2,\ldots \in V$ with $\lim u_n = u$ weakly in $V$.

Clearly, if the sequence $(\gota_n)_{n \in \Ni}$ is uniformly $\tilde j$-elliptic
and converges weakly to $\gota$, then the form $\gota$ is also $\tilde j$-elliptic
and satisfies the bounds (\ref{eSdtonl5;1}) with $\gota_n$ replaced by~$\gota$.
Then $A_n$ and $A$ are self-adjoint graphs for all $n \in \Ni$ by Theorem~\ref{tdtonl202}.
A natural question is whether these conditions suffice to show 
$\lim_{n \to \infty} (A_n + i \, s \, I)^{-1} = (A + i \, s \, I)^{-1}$
strongly for all $s \in \Ri \setminus \{ 0 \} $.

There is a surprisingly simple counter-example which shows that 
more conditions are needed.

\begin{exam} \label{xdton501}
Choose $V = \widetilde H = \Ci^2$, $H = \Ci$, 
$j(u) = u_1$, $\tilde j(u) = u$, 
$\gota(u,v) = 0$ and 
\[
\gota_n(u,v)
= \tfrac{1}{n} (u_1 \, \overline {v_2} + u_2 \, \overline{v_1}) 
\]
for all $n \in \Ni$.
Clearly $\tilde j$ is compact and the sequence $(\gota_n)_{n \in \Ni}$ is 
uniformly $\tilde j$-elliptic and converges weakly to $\gota$.
An easy calculation gives 
$A = \Ci \times \{ 0 \} $ and $A_n = \{ 0 \} \times \Ci$
for all $n \in \Ni$.
Note that $A_n$ is multi-valued.
If $s \in \Ri \setminus \{ 0 \} $ then
\[
(A + i \, s \, I)^{-1} = \tfrac{1}{i \, s} \, I
\quad \mbox{and} \quad
(A_n + i \, s \, I)^{-1} = 0
\]
for all $n \in \Ni$.
Therefore $\lim_{n \to \infty} (A_n + i \, s \, I)^{-1} \neq (A + i \, s \, I)^{-1}$
in any Hausdorff topology on $\cl(H)$.
\end{exam}

To understand this counter-example better, we recall the spaces of non-uniqueness
\begin{eqnarray*}
W(\gota_n) & = & \{ u \in \ker j : \gota_n(u,v) = 0 \mbox{ for all } v \in V \} \mbox{ and}  \\
W(\gota) & = & \{ u \in \ker j : \gota(u,v) = 0 \mbox{ for all } v \in V \}
.
\end{eqnarray*}
In Example~\ref{xdton501} one has $\dim W(\gota) = 1$ whilst $\dim W(\gota_n) = 0$ for all $n \in \Ni$.
We shall show in Proposition~\ref{pdtonl506} that 
in general $\limsup_{n \to \infty} \dim W(\gota_n) \leq \dim W(\gota)$
if the sequence $(\gota_n)_{n \in \Ni}$ is uniformly $\tilde j$-elliptic and converges 
weakly to $\gota$.
The first main theorem of this section is the following.

\begin{thm} \label{tdtonl502.3}
For all $n \in \Ni$ let $\gota, \gota_n \colon V \times V \to \Ci$ be continuous symmetric forms.
Suppose that the sequence $(\gota_n)_{n \in \Ni}$ is uniformly $\tilde j$-elliptic and converges 
weakly to $\gota$.
Moreover, suppose that $W(\gota) = \{ 0 \} $.
Then 
\[
\lim_{n \to \infty} (A_n + i \, s \, I)^{-1} 
= (A + i \, s \, I)^{-1}
\]
strongly for all $s \in \Ri \setminus \{ 0 \} $.
Moreover, if in addition the map $j$ is compact, then the convergence is uniform
in $\cl(H)$.
\end{thm}

This theorem will be a special case of Theorem~\ref{tdtonl502} together with 
Proposition~\ref{pdtonl506}, which we prove later.
It will allow us to prove convergence results, not only for 
our basic example,
but also for Dirichlet-to-Neumann graphs associated with elliptic operators
(see Section~\ref{Sdtonl7}).

In Theorem~\ref{tdtonl502.3} one has uniform resolvent convergence if $j$ is compact. 
Moreover, each graph $A_n$ is lower bounded by Theorem~\ref{tdtonl210.5}.
Hence it is tempting to conjecture that the graphs $A_n$ are 
lower bounded uniformly in $n \in \Ni$.
The next example shows that this conjecture is false in general.

\begin{exam} \label{xdtonl508}
Choose $V = \widetilde H = \Ci^2$, $H = \Ci$, 
$j(u) = u_1$, $\tilde j(u) = u$, 
$\gota(u,v) = u_1 \, \overline {v_2} + u_2 \, \overline{v_1}$ and 
\[
\gota_n(u,v)
= u_1 \, \overline {v_2} + u_2 \, \overline{v_1} + \tfrac{1}{n} \, u_2 \, \overline{v_2}
\]
for all $n \in \Ni$.
Clearly $\tilde j$ is compact and the sequence $(\gota_n)_{n \in \Ni}$ is 
uniformly $\tilde j$-elliptic and converges weakly to $\gota$.
Moreover, $\ker j = \{ 0 \} \times \Ci$
and $W(\gota) = W(\gota_n) = \{ 0 \} $ for all $n \in \Ni$.
So by Theorem~\ref{tdtonl502.3} the sequence $(A_n)$ converges 
in the uniform resolvent sense to $A$.
But 
\[
A_n = \{ (\lambda, -n \lambda) : \lambda \in \Ci \} 
\]
for all $n \in \Ni$.
Hence $A_n$ is not lower bounded uniformly in $n \in \Ni$.
\end{exam}

In Example~\ref{xdtonl508} one has 
\begin{eqnarray*}
V(\gota) & = & \{ (0,\lambda) : \lambda \in \Ci \} \mbox{ and}  \\
V(\gota_n) & = & \{ (-\frac{1}{n} \, \lambda, \lambda) : \lambda \in \Ci \} 
\end{eqnarray*}
for all $n \in \Ni$.
So $V(\gota_n) \cap \ker j = \{ 0 \} $ and $V(\gota) \cap \ker j = \{ 0 \} \times \Ci$.
Hence $\dim (V(\gota_n) \cap \ker j) = 0$ for all $n \in \Ni$ whilst 
$\dim (V(\gota) \cap \ker j) = 1$.
This, together with Example~\ref{xdton501} and Theorem~\ref{tdtonl502.3}, 
suggests that the dimension of the spaces $W(\gota)$ and $W(\gota_n)$ is intimately 
connected with the question of whether the $A_n$ converge in the strong resolvent 
sense, while the dimension of $V(\gota) \cap \ker j$ and $V(\gota_n) \cap \ker j$ 
influences uniform lower boundedness.
This will be the subject of Theorems~\ref{tdtonl502} and \ref{tdtonl509}, respectively.
As the second main theorem of this section, we first give a special case of 
Theorem~\ref{tdtonl509}, which is analogous to Theorem~\ref{tdtonl502.3}.

\begin{thm} \label{tdtonl509.5}
For all $n \in \Ni$ let $\gota, \gota_n \colon V \times V \to \Ci$ be 
continuous symmetric forms.
Suppose that the sequence $(\gota_n)_{n \in \Ni}$ is uniformly $\tilde j$-elliptic and converges 
weakly to $\gota$.
Moreover, suppose that $V(\gota) \cap \ker j = \{ 0 \} $.
Then the graphs $A_n$ are bounded below uniformly in $n \in \Ni$.
\end{thm}

We now wish to develop the prerequisites necessary for the proofs of 
Theorems~\ref{tdtonl502} and \ref{tdtonl509}. These are quite similar 
and much of what follows will be used for both.
Throughout the remainder of this section we assume,
in addition to the assumption that $\gota_n$ and $\gota$ are continuous and 
symmetric for all $n \in \Ni$, 
that the sequence $(\gota_n)_{n \in \Ni}$ is uniformly $\tilde j$-elliptic and converges 
weakly to $\gota$.

The first lemma is of independent interest.

\begin{lemma} \label{ldtonl503.1}
Let $(u_n)_{n \in \Ni}$ be a sequence in $V$ and $u \in V$.
Suppose that $\lim_{n \to \infty} u_n = u$ weakly in $V$ and 
$\lim_{n \to \infty} \gota_n(u_n) = \gota(u)$.
Then $\lim_{n \to \infty} u_n = u$ strongly in $V$.
\end{lemma}
\proof\
Let $n \in \Ni$.
Then 
\[
\gota_n(u_n - u) 
= \gota_n(u_n) - 2 \RRe \gota_n(u_n,u) + \gota_n(u)
.  \]
So $\lim \gota_n(u_n - u) = 0$ by assumption and the weak convergence of $(\gota_n)_{n \in \Ni}$.
Clearly $\lim \tilde j(u_n - u) = 0$ in $\widetilde H$.
Finally, the uniform $\tilde j$-ellipticity (\ref{eSdtonl5;1}) gives
\[
\mu \, \|u_n - u\|_V^2 
\leq \gota_n(u_n - u) + \omega \, \|\tilde j(u_n - u)\|_{\widetilde H}^2
\]
for all $n \in \Ni$ and the lemma follows.\hfill$\Box$

\begin{lemma} \label{ldtonl503}
Suppose that either 
\begin{eqnarray*}
&  \quad & U = W(\gota) \mbox{ and } U_n = W(\gota_n) \mbox{ for all } n \in \Ni \mbox{, or,}  \\
&  \quad & U = V(\gota) \cap \ker j \mbox{ and } U_n = V(\gota_n) \cap \ker j \mbox{ for all } n \in \Ni .
\end{eqnarray*}
For all $n \in \Ni$ let $u_n \in U_n$ and let $u \in V$.
Suppose that $\lim_{n \to \infty} u_n = u$ weakly in $V$.
Then $u \in U$ and $\lim_{n \to \infty} u_n = u$ strongly in $V$.
\end{lemma}
\proof\
Since $j$ is weakly continuous and $u_n \in U_n \subset \ker j$ for all $n \in \Ni$,
one deduces that $j(u) = \lim j(u_n) = 0$.
Moreover, $\gota(u,v) = \lim \gota_n(u_n,v)$ for all $v \in V$.
So $u \in U$ in both cases.
In particular, $\gota(u) = 0$.
Clearly $\gota_n(u_n) = 0$ for all $n \in \Ni$.
Now use Lemma~\ref{ldtonl503.1}.\hfill$\Box$

\vertspace

By Remark~\ref{rdtonl202.7} we know that the spaces $W(\gota)$ and $W(\gota_n)$ are finite dimensional
for all $n \in \Ni$.
The same argument also shows that the spaces $V(\gota) \cap \ker j$ and 
$V(\gota_n) \cap \ker j$ are finite dimensional.
The weak convergence of $(\gota_n)$ allows one to compare their dimensions in the next proposition.

\begin{prop} \label{pdtonl505}
Suppose that either 
\begin{eqnarray*}
&  \quad & U = W(\gota) \mbox{ and } U_n = W(\gota_n) \mbox{ for all } n \in \Ni \mbox{, or,}  \\
&  \quad & U = V(\gota) \cap \ker j \mbox{ and } U_n = V(\gota_n) \cap \ker j \mbox{ for all } n \in \Ni .
\end{eqnarray*}
Then  
\begin{equation}
\limsup_{n \to \infty} \dim U_n \leq \dim U
.
\label{epdtonl505;10}
\end{equation}
Moreover, if  $\dim U_n = \dim U$ for all $n \in \Ni$, then for all $u \in U$ 
there exists a 
sequence $(u_n)_{n \in \Ni}$ in $V$ such that $u_n \in U_n$ for all $n \in \Ni$ 
and $\lim_{n \to \infty} u_n = u$.
\end{prop}
\proof\
Let $d_0 \in \Ni$ and suppose that $d_0 \leq \limsup \dim U_n$.
We shall prove that $d_0 \leq \dim U$.
This implies (\ref{epdtonl505;10}).

Without loss of generality we may assume that $d_0 \leq \dim U_n$ for all $n \in \Ni$.
For all $n \in \Ni$ let $ \{ u_{n1},\ldots,u_{n d_0} \} $ be an orthonormal 
set in $U_n$ of dimension $d_0$.
Let $\ell \in \{ 1,\ldots,d_0 \} $.
Then $(u_{n \ell})_{n \in \Ni}$ is a bounded sequence in $V$, so 
passing to a subsequence if necessary, there exists a $u_\ell \in V$
such that $\lim u_{n\ell} = u_\ell$ weakly in $V$.
Then $u_\ell \in U$ and $\lim u_{n\ell} = u_\ell$ strongly in $V$ by 
Lemma~\ref{ldtonl503}.
Since $ \{ u_{n1},\ldots,u_{n d_0} \} $ is an orthonormal 
set in $V$ for all $n \in \Ni$, also 
$ \{ u_1,\ldots,u_{d_0} \} $ is an orthonormal set in $V$, of dimension $d_0$.
Hence $\dim U \geq d_0$.

For the last statement choose $d_0 = \dim U$ and fix $u \in U$.
Then the above gives that there exists a subsequence $(U_{n_k})_{k \in \Ni}$ 
of $(U_n)_{n \in \Ni}$ and for all $k \in \Ni$ there exists a $u_k \in U_{n_k}$
such that $\lim_{k \to \infty} u_k = u$.
Hence for all $\varepsilon > 0$ there exists an $N \in \Ni$ such that 
for all $n \in \Ni$ with $n \geq N$ there is a $v \in U_n$ such that 
$\|v - u\| < \varepsilon$.
This implies the last statement.\hfill$\Box$

\vertspace

Proposition~\ref{pdtonl505} gives a remarkable inequality for the dimensions 
of the spaces of non-uniqueness if the sequence of forms $(\gota_n)$ converges weakly.
We state it explicitly.

\begin{prop} \label{pdtonl506}
 $\limsup\limits_{n \to \infty} \dim W(\gota_n) \leq \dim W(\gota)$ and also
$\limsup\limits_{n \to \infty} \dim (V(\gota_n) \cap \ker j) \leq \dim (V(\gota) \cap \ker j)$.
\end{prop}

\begin{prop} \label{pdtonl507}
Suppose that either 
\begin{eqnarray*}
&  \quad & U = W(\gota) \mbox{ and } U_n = W(\gota_n) \mbox{ for all } n \in \Ni \mbox{, or}  \\
&  \quad & U = V(\gota) \cap \ker j \mbox{ and } U_n = V(\gota_n) \cap \ker j \mbox{ for all } n \in \Ni ,
\end{eqnarray*}
and that $\lim_{n \to \infty} \dim U_n = \dim U$.
Let $u \in V$ and for all $n \in \Ni$ let $u_n \in U_n^\perp$.
If $\lim_{n \to \infty} u_n = u$ weakly in $V$, then $u \in U^\perp$.
\end{prop}
\proof\
Let $v \in U$.
By Proposition~\ref{pdtonl505}
for all $n \in \Ni$ there exists a $v_n \in U_n$ such that 
$\lim v_n = v$ strongly in $V$.
Then $(u_n,v_n)_V = 0$ for all $n \in \Ni$.
Taking the limit $n \to \infty$ one deduces that 
$(u,v)_V = \lim (u_n,v_n)_V = 0$.
So $u \in U^\perp$.\hspace*{5mm}\hfill$\Box$

\vertspace

Now we are able to prove our main convergence result, which is 
the extension of Theorem~\ref{tdtonl502.3} to which we have alluded.

\begin{thm} \label{tdtonl502}
For all $n \in \Ni$ let $\gota, \gota_n \colon V \times V \to \Ci$ be continuous symmetric forms.
Suppose that the sequence $(\gota_n)_{n \in \Ni}$ is uniformly $\tilde j$-elliptic and converges 
weakly to $\gota$.
Moreover, suppose that $\lim_{n \to \infty} \dim W(\gota_n) = \dim W(\gota)$.
Then 
\[
\lim_{n \to \infty} (A_n + i \, s \, I)^{-1} 
= (A + i \, s \, I)^{-1}
\]
strongly for all $s \in \Ri \setminus \{ 0 \} $.
Moreover, if in addition the map $j$ is compact, then the convergence is uniform
in $\cl(H)$.
\end{thm}
\proof\
Let $y,y_1,y_2,\ldots \in H$ and suppose that $\lim y_n = y$ weakly in $H$.
For all $n \in \Ni$ define $x_n = (A_n + i \, s \, I)^{-1} y_n$.
There exists a $u_n \in V$ such that $j(u_n) = x_n$ and 
\begin{equation}
\gota_n(u_n,v) = (y_n - i \, s \, x_n,j(v))_H
\label{etdtonl213;2}
\end{equation}
for all $v \in V$.
Without loss of generality we may assume that $u_n \in W(\gota_n)^\perp$
by (\ref{etdton202;1}).
Then $\|x_n\|_H \leq \frac{1}{|s|} \, \|y_n\|_H$ and 
$|\gota_n(u_n)| = |(y_n - i \, s \, x_n, j(u_n))_H| \leq \frac{2}{|s|} \, \|y_n\|_H^2$.
Since the sequence $(y_n)_{n \in \Ni}$ converges weakly, it is bounded.
Let $M > 0$ be such that $\|y_n\|_H \leq M$ for all $n \in \Ni$.
We shall prove that the sequence $(u_n)_{n \in \Ni}$ is bounded in $V$.
If not, then it follows from (\ref{eSdtonl5;1}) that 
the sequence $(\tilde j(u_n))_{n \in \Ni}$ is not bounded in 
$\widetilde H$.
Passing to a subsequence if necessary, we may assume that 
$\lim \|\tilde j(u_n)\|_{\widetilde H} = \infty$.
Write $\lambda_n = \|\tilde j(u_n)\|_{\widetilde H}$ for all $n \in \Ni$.
Using again (\ref{eSdtonl5;1}) it follows that 
$(\frac{1}{\lambda_n} \, u_n)_{n \in \Ni}$ is bounded in $V$.
Passing to a subsequence if necessary, it follows that there exists 
a $u_0 \in V$ such that $\lim \frac{1}{\lambda_n} \, u_n = u_0$ weakly in $V$.
Then $\lim \frac{1}{\lambda_n} \, \tilde j(u_n) = \tilde j(u_0)$ strongly 
in $\widetilde H$, since $\tilde j$ is compact.
Moreover, $\|\tilde j(u_0)\|_{\widetilde H} = 1$.
Using (\ref{eSdtonl5;2}) and (\ref{etdtonl213;2}) one deduces that 
\[
\gota(u_0,v)
= \lim_{n \to \infty} \gota_n(\tfrac{1}{\lambda_n} \, u_n, v)
= \lim_{n \to \infty} \tfrac{1}{\lambda_n} \, (y_n - i \, s \, x_n,j(v))_H
= 0
\]
for all $v \in V$.
Moreover, $j(u_0) = \lim j(\tfrac{1}{\lambda_n} \, u_n)
= \lim \tfrac{1}{\lambda_n} \, x_n = 0$.
So $u_0 \in W(\gota)$.
But also $u_0 \in W(\gota)^\perp$ by Proposition~\ref{pdtonl507}.
Hence $u_0 = 0$.
But $\|\tilde j(u_0)\|_{\widetilde H} = 1$.
This is a contradiction.
Hence the sequence $(u_n)_{n \in \Ni}$ is bounded in $V$.

Passing to a subsequence if necessary, there exists a $u_0 \in V$ 
such that $\lim u_n = u_0$ weakly in $V$.
Set $x_0 = j(u_0)$.
Then $\lim x_n = \lim j(u_n) = j(u_0) = x_0$ weakly in $H$.
Let $v \in V$.
Using (\ref{eSdtonl5;2}) and (\ref{etdtonl213;2}) one deduces that 
\begin{equation}
\gota(u_0,v)
= \lim_{n \to \infty} \gota_n(u_n,v)
= \lim_{n \to \infty} (y_n - i \, s \, x_n, j(v))_H
= (y - i \, s \, x_0, j(v))_H
\label{etdtonl213;5} 
\end{equation}
for all $v \in V$.
So $(x_0, y - i \, s \, x_0) \in A$.
Then $(A + i \, s \, I)^{-1} y = x_0$.

Now we prove the part of the theorem concerning strong convergence.
Let $y \in H$.
Choose $y_n = y$ for all $n \in \Ni$.
It follows from (\ref{etdtonl213;2}) that 
$\gota_n(u_n) = (y - i \, s \, x_n,x_n)_H = (y,x_n)_H - i \, s \, \|x_n\|_H^2$ for all $n \in \Ni$.
Since $\gota$ is symmetric, one has $\gota_n(u_n) \in \Ri$.
Therefore $s \, \|x_n\|_H^2 = \IIm (y,x_n)_H$ for all $n \in \Ni$.
Similarly, $s \, \|x_0\|_H^2 = \IIm (y,x_0)_H$ by (\ref{etdtonl213;5}).
Since $\lim (y,x_n)_H = (y,x_0)_H$ by the weak convergence of 
$(x_n)_{n \in \Ni}$, one has 
$\lim \|x_n\|_H^2 = \|x_0\|_H^2$.
Hence $\lim x_n = x_0$ strongly in $H$.
This implies the strong resolvent convergence.

Finally suppose that the map $j$ is compact.
We shall prove that 
\[
\lim_{n \to \infty} (A_n + i \, s \, I)^{-1} 
= (A + i \, s \, I)^{-1} 
\]
in $\cl(H)$.
Suppose not. 
Then passing to a subsequence if necessary,
there are $\varepsilon > 0$ and $y_1,y_2,\ldots \in H$ 
such that $\|y_n\|_H \leq 1$ and 
$\|(A_n + i \, s \, I)^{-1} y_n - (A + i \, s \, I)^{-1} y_n\|_H \geq \varepsilon$
for all $n \in \Ni$.
Passing to a subsequence if necessary, there exists a $y \in H$ such that 
$\lim y_n = y$ weakly in $H$.
For all $n \in \Ni$ there exists a $u_n \in V$ such that 
$j(u_n) = (A_n + i \, s \, I)^{-1} y_n$ and 
$\gota_n(u_n,v) = (y_n - i \, s \, j(u_n), j(v))_H$ for all $v \in V$.
By the above there exists a $u \in V$ such that,
again passing to a subsequence if necessary, 
$\lim u_n = u$ weakly in $V$ and 
$j(u) = (A + i \, s \, I)^{-1} y$.
Since $j$ is compact one has 
\[
\lim (A_n + i \, s \, I)^{-1} y_n
= \lim j(u_n)
= j(u)
= (A + i \, s \, I)^{-1} y
.  \]
Moreover, the operator $(A + i \, s \, I)^{-1}$ is compact by Proposition~\ref{tdtonl202}.
Therefore 
\[
\lim (A + i \, s \, I)^{-1} y_n = (A + i \, s \, I)^{-1} y
.  \]
This is a contradiction.
The proof of Theorem~\ref{tdtonl502} is complete.\hfill$\Box$

\begin{remarkn} \label{rdtonl502.1}
Note that actually $\lim u_n = u_0$ strongly in $V$ in the proof of 
Theorem~\ref{tdtonl502}.
The reason is as follows.
If $n \in \Ni$ then 
$\gota_n(u_n) + i \, s \, \|x_n\|_H^2 = (y,x_n)_H$ and 
$\gota(u_0) + i \, s \, \|x_0\|_H^2 = (y,x_0)_H$.
So $\lim \gota_n(u_n) = \gota(u_0)$.
Now apply Lemma~\ref{ldtonl503.1}.
\end{remarkn}

We now prove our main result on uniform lower boundedness, the 
analogue of Theorem~\ref{tdtonl502} involving the spaces 
$V(\gota_n) \cap \ker j$ rather than $W(\gota_n)$.
It turns out that this property will be needed to prove convergence 
of the associated semigroups in Section~\ref{tdtonl6}.

\begin{thm} \label{tdtonl509}
For all $n \in \Ni$ let $\gota,\gota_n \colon V \times V \to \Ci$ be
continuous symmetric forms.
Suppose that the sequence $(\gota_n)_{n \in \Ni}$ is uniformly $\tilde j$-elliptic and converges 
weakly to $\gota$.
Moreover, suppose that $\lim_{n \to \infty} \dim (V(\gota_n) \cap \ker j) = \dim (V(\gota) \cap \ker j)$.
Then the graphs $A_n$ are bounded below uniformly in $n \in \Ni$.
\end{thm}
\proof\
Suppose that the graphs $A_n$ are not bounded below uniformly in $n \in \Ni$.
Then for all $n \in \Ni$ there exists a pair $(x_n,y_n) \in A_n$ such that 
$(y_n , x_n)_H + n \, \|x_n\|_H^2 < 0$.
By Corollary~\ref{cdtonl207} for all $n \in \Ni$ there exists a 
$u_n \in V(\gota_n) \cap (V(\gota_n) \cap \ker j)^\perp$ such that 
$j(u_n) = x_n$ and 
$\gota_n(u_n,v) = (y_n,j(v))_H$ for all $v \in V(\gota_n) \cap (V(\gota_n) \cap \ker j)^\perp$.
Then $u_n \neq 0$.
Without loss of generality we may assume that $\|\tilde j(u_n)\|_{\widetilde H} = 1$.
Let $\mu,\omega > 0$ be as in (\ref{eSdtonl5;1}).
Then 
\[
\mu \, \|u_n\|_V^2 + n \, \|j(u_n)\|_H^2
\leq \gota_n(u_n) + \omega \, \|\tilde j(u_n)\|_{\widetilde H}^2 + n \, \|j(u_n)\|_H^2  
< \omega \, \|\tilde j(u_n)\|_{\widetilde H}^2
= \omega
\]
for all $n \in \Ni$.
Hence the sequence $(u_n)_{n \in \Ni}$ is bounded in $V$ and $\lim j(u_n) = 0$ in $H$.
Passing to a subsequence if necessary we may assume that there exists a $u \in V$
such that $\lim u_n = u$ weakly in $V$.
Then $\|\tilde j(u)\|_{\widetilde H} = 1$ since $\tilde j$ is compact.
Moreover, $j(u) = 0$, so $u \in \ker j$.

Let $v \in \ker j$. 
If $n \in \Ni$, then $\gota_n(u_n,v) = 0$ since $u_n \in V(\gota_n)$.
Hence $\gota(u,v) = \lim \gota_n(u_n,v) = 0$.
So $u \in V(\gota)$.
Thus $u \in V(\gota) \cap \ker j$.

Next, $u_n \in (V(\gota_n) \cap \ker j)^\perp$ for all $n \in \Ni$ and 
$\lim u_n = u$ weakly in $V$. 
Therefore Proposition~\ref{pdtonl507} implies that $u \in (V(\gota) \cap \ker j)^\perp$.
Hence $u = 0$.
But $\|\tilde j(u)\|_{\widetilde H} = 1$.
This is a contradiction.
Hence the graphs $A_n$ are bounded below uniformly in $n \in \Ni$.\hfill$\Box$

\vertspace

Now Theorem~\ref{tdtonl509.5} is an easy corollary.
Note that Theorem~\ref{tdtonl509} also gives a different proof of Theorem~\ref{tdtonl210.5}.

We next wish to compare briefly the conditions on the dimensions
in Theorems~\ref{tdtonl502} and \ref{tdtonl509} through examples, 
in particular as regards sufficiency and necessity.
We first observe that the conditions in Theorem~\ref{tdtonl502}, while 
sufficient (and arguably somehow natural) for resolvent convergence,
are not necessary, as the following example shows.

\begin{exam} \label{xdtonl512}
Let $V = \Ci^2$, $H = \Ci$ and define $j \colon V \to H$ by $j(u) = u_1$.
Define the forms $\gota$ and $\gota_n$ on $V$ by 
$\gota(u,v) = u_1 \, \overline{v_1}$ and 
$\gota_n(u,v) = u_1 \, \overline{v_1} + \frac{1}{n} \, u_2 \, \overline{v_2}$ for all $n \in \Ni$.
Then $A = A_n = I$, $W(\gota) = \{ 0 \} \times \Ci$ and $W(\gota_n) = \{ 0 \} $
for all $n \in \Ni$.
So clearly $\lim A_n = A$ uniformly in resolvent sense, but 
$\lim \dim W(\gota_n) < \dim W(\gota)$.
\end{exam}

Similarly, the conditions in Theorem~\ref{tdtonl509} are sufficient for 
uniform lower boundedness, but not necessary.
An example is as follows (cf.\ Example~\ref{xdtonl508}).

\begin{exam} \label{xdtonl512.4}
Choose $V = \widetilde H = \Ci^2$, $H = \Ci$, 
$j(u) = u_1$, $\tilde j(u) = u$, 
$\gota(u,v) = u_1 \, \overline {v_2} + u_2 \, \overline{v_1}$ and 
\[
\gota_n(u,v)
= u_1 \, \overline {v_2} + u_2 \, \overline{v_1} - \tfrac{1}{n} \, u_2 \, \overline{v_2}
\]
for all $n \in \Ni$.
Clearly $\tilde j$ is compact and the sequence $(\gota_n)_{n \in \Ni}$ is 
uniformly $\tilde j$-elliptic and converges weakly to $\gota$.
Moreover, $V(\gota_n) \cap \ker j = \{ 0 \} $ for all $n \in \Ni$
and $V(\gota) \cap \ker j = \{ 0 \} \times \Ci$.
Hence $\dim (V(\gota_n) \cap \ker j) = 0$ for all $n \in \Ni$ whilst 
$\dim (V(\gota) \cap \ker j) = 1$.
But
\[
A_n = \{ (\lambda, n \lambda) : \lambda \in \Ci \} 
\]
for all $n \in \Ni$.
Hence the $A_n$ are lower bounded uniformly in $n \in \Ni$.
\end{exam}

We also note that the two conditions are not interchangeable.
The following counter-example shows that
$\lim_{n \to \infty} \dim (V(\gota_n) \cap \ker j) = \dim (V(\gota) \cap \ker j)$
is not sufficient for strong resolvent converge of the sequence $(A_n)$ to $A$.
We recall meanwhile that Example~\ref{xdtonl508} shows that the condition
$\lim_{n \to \infty} \dim W(\gota_n) = \dim W(\gota)$
is likewise not sufficient for uniform lower boundedness.

\begin{exam} \label{xdtonl5010}
Let  $V$, $H$, $\widetilde H$, $j$, $\tilde j$, $\gota$ and $\gota_n$ 
be as in Example~\ref{xdton501}.
Then the sequence of graphs $(A_n)$ does not converge to $A$ in the resolvent sense.
But $\ker j = \{ 0 \} \times \Ci$, $V(\gota_n) = \ker j$ and $V(\gota) = V$ for 
all $n \in \Ni$.
So $\dim (V(\gota_n) \cap \ker j) = 1 = \dim (V(\gota) \cap \ker j)$ for all $n \in \Ni$
and the graphs $A_n$ are bounded below uniformly in $n \in \Ni$.
\end{exam}

It is possible to reformulate the conditions on the dimensions
$\lim_{n \to \infty} \dim W(\gota_n) = \dim W(\gota)$
and 
$\lim_{n \to \infty} \dim (V(\gota_n) \cap \ker j) = \dim (V(\gota) \cap \ker j)$
using the concept of the gap between two closed subspaces.
We define the gap $\hat \delta(M,N)$ between two closed subspaces
$M$ and $N$ of $H$ as in Kato \cite{Kat1} (IV.2.2),
that is,
\begin{eqnarray*}
\delta(M,N) & = & \sup_{\scriptstyle u \in M \atop
                        \scriptstyle \|u\| \leq 1}
                     d(u,N),  \\
\hat \delta(M,N) & = & \max( \delta(M,N), \delta(N,M)).
\end{eqnarray*}

\begin{lemma} \label{ldtonl511}
Suppose that either
\begin{eqnarray*}
&  \quad & U = W(\gota) \mbox{ and } U_n = W(\gota_n) \mbox{ for all } n \in \Ni \mbox{, or}  \\
&  \quad & U = V(\gota) \cap \ker j \mbox{ and } U_n = V(\gota_n) \cap \ker j \mbox{ for all } n \in \Ni .
\end{eqnarray*}
Then $\lim_{n \to \infty} \delta(U_n,U) = 0$.
Moreover, let $P_U$ and $P_{U_n}$ be the orthogonal projection in $H$ onto $U$ and $U_n$
for all $n \in \Ni$.
Then the following are equivalent.
\begin{tabeleq}
\item \label{ldtonl511-1}
$\lim\limits_{n \to \infty} \dim U_n = \dim U$.
\item \label{ldtonl511-2}
$\lim\limits_{n \to \infty} \hat\delta(U_n,U) = 0$.
\item \label{ldtonl511-3}
$\lim\limits_{n \to \infty} \delta(U,U_n) = 0$.
\item \label{ldtonl511-4}
$\lim\limits_{n \to \infty} P_{U_n} = P_U$ in $\cl(H)$.
\end{tabeleq}
\end{lemma}
\proof\
We first show that $\lim_{n \to \infty} \delta(U_n,U) = 0$.
Suppose not. 
Then passing to a subsequence if necessary, there exists an $\varepsilon > 0$ 
such that $\delta(U_n,U) \geq 2\varepsilon$ for all $n \in \Ni$.
Then for all $n \in \Ni$ there exists a $u_n \in U_n$ such that 
$d(u_n,U) \geq \varepsilon$ and $\|u_n\| \leq 1$.
Passing to a subsequence if necessary, there exists a $u \in V$ 
such that $\lim u_n = u$ weakly in $V$.
By Lemma~\ref{ldtonl503} it follows that $u \in U$ and $\lim u_n = u$
strongly in $V$.
So $\varepsilon \leq d(u_n,U) \leq \|u_n - u\|$ for all $n \in \Ni$.
This is a contradiction.

The equivalence `\ref{ldtonl511-2}$\Leftrightarrow$\ref{ldtonl511-3}' is now 
trivial.

`\ref{ldtonl511-1}$\Rightarrow$\ref{ldtonl511-3}'.
Suppose \ref{ldtonl511-1} and not \ref{ldtonl511-3}. 
Then passing to a subsequence if necessary, there exists an $\varepsilon > 0$ 
and for all $n \in \Ni$ there exists a $v_n \in U$ such that 
$\|v_n\| \leq 1$ and $d(v_n,U_n) \geq \varepsilon$.
The space $U$ is finite dimensional.
So passing to a subsequence if necessary, there exists a $u \in U$ 
such that $\lim v_n = u$ strongly in $U$.
By Proposition~\ref{pdtonl505} and \ref{ldtonl511-1},
for all $n \in \Ni$ there exists a 
$u_n \in U_n$ such that $\lim u_n = u$ in $V$.
Then 
$
d(v_n,U_n)
\leq \|v_n - u_n\|
\leq \|v_n - u\| + \|u - u_n\|$ for all $n \in \Ni$.
So $\lim d(v_n,U_n) = 0$. 
This is a contradiction.

`\ref{ldtonl511-2}$\Rightarrow$\ref{ldtonl511-1}'.
This follows from \cite{Kat1} Corollary~IV.2.6.

`\ref{ldtonl511-2}$\Leftrightarrow$\ref{ldtonl511-4}'.
This is in \cite{Kat1} Footnote~1 on page~198.\hfill$\Box$

\section{Semigroup convergence} \label{tdtonl6}

In this section we consider the semigroup generated by a self-adjoint graph
which is bounded below.
Let $A$ be a self-adjoint graph which is bounded below.
For all $t > 0$ define the operator $e^{-tA}$ by
\[
e^{-tA} = \lim_{n \to \infty} \Big( (I + \tfrac{t}{n} \, A)^{-1} \Big)^n
.  \]
Using the decomposition $H = A(0) \oplus A(0)^\perp$ 
one has $e^{-tA} = 0 \oplus e^{-t A^\circ}$.
We call $(e^{-tA})_{t > 0}$ the {\bf semigroup generated by $A$}.

Let $A,A_1,A_2,\ldots$ be self-adjoint graphs.
In Theorem~\ref{tdtonl502} we provided conditions such that 
$(A_n)$ converges uniformly to $A$ in the resolvent sense.
So one might hope that then the semigroups converge too, at least pointwise.
In general, however, this is false.

\begin{exam} \label{xdtonl603}
Consider Example~\ref{xdtonl508} again.
Then $H = \Ci$ and $e^{-t A_n} 1 = e^{nt}$ for all $t > 0$ and $n \in \Ni$.
So $\lim e^{-t A_n} x$ does not exist in $H$ for any $x \in H \setminus \{ 0 \} $ and 
any $t > 0$.
\end{exam}

The main problem in Example~\ref{xdtonl603} is that the sequence of
self-adjoint graphs is not bounded below uniformly in $n \in \Ni$.

\begin{lemma} \label{ldtonl604}
Let $A,A_1,A_2,\ldots$ be positive self-adjoint graphs.
Suppose that $\lim_{n \to \infty} A_n = A$ in the strong resolvent sense.
Let $K \subset \Ci \setminus [0,\infty)$ be compact and $y \in H$.
Then 
\[
\lim_{n \to \infty} (A_n + \lambda \, I)^{-1} y = (A + \lambda \, I)^{-1} y
\]
in $H$ uniformly for all $\lambda \in K$.
\end{lemma} 
\proof\
Let $n \in \Ni$. 
Then $\|(A_n + \lambda \, I)^{-1}\| \leq \frac{1}{|\IIm \lambda|}$
if $\lambda \in \Ci \setminus \Ri$ and 
$\|(A_n + \lambda \, I)^{-1}\| \leq \frac{1}{\lambda}$
if $\lambda \in (0,\infty)$.
Similar estimates are valid for $A$.
The resolvent identities (\ref{eSdtonl3.1;1}) then imply 
locally uniform convergence in $\lambda$,
and the lemma now follows by a compactness argument.\hspace*{5mm}\hfill$\Box$

\begin{thm} \label{tdtonl605}
Let $A,A_1,A_2,\ldots$ be self-adjoint graphs which are uniformly bounded below.
Suppose that $\lim_{n \to \infty} A_n = A$ in the strong resolvent sense.
Then $\lim_{n \to \infty} e^{-t A_n} = e^{-t A}$ strongly for all $t > 0$.
An analogous statement is valid for uniform convergence.
\end{thm}
\proof\
We may assume that the $A_n$ and $A$ are positive.
Let $\gamma$ be the contour in $\Ci$ formed by combining 
the two line segments 
$ \{ \lambda \in \Ci : \arg \lambda = \pm \frac{3\pi}{2} \mbox{ and } |\lambda| \geq 1 \} $
together with the arc
$ \{ \lambda \in \Ci : |\lambda| = 1 \mbox{ and } |\arg \lambda| \leq \frac{3\pi}{2} \} $.
Let $n \in \Ni$.
Then $\|(A_n + \lambda \, I)^{-1}\| \leq \frac{1}{|\IIm \lambda|}$
if $\lambda \in \Ci \setminus \Ri$ and a similar estimate is valid for $A$.
So 
\[
e^{-t A_n} y = \frac{1}{2 \pi i } \int_\gamma e^{t \lambda} (A_n + \lambda \, I)^{-1} y \, d\lambda
\]
for all $t > 0$ and $y \in H$.
Now take the limit $n \to \infty$ and use Lemma~\ref{ldtonl604}.\hfill$\Box$

\section{Convergence of Dirichlet-to-Neumann resolvents \\and semigroups} \label{Sdtonl7}

In this section we give applications of the generation result,
Theorem~\ref{tdtonl202}, and the convergence theorem, 
Theorem~\ref{tdtonl502.3}.
Throughout this section $\Omega \subset \Ri^d$ is Lipschitz
with boundary $\Gamma = \partial \Omega$ and we choose 
$j \colon H^1(\Omega) \to L_2(\Gamma)$ the trace operator.
Moreover, $\tilde j \colon H^1(\Omega) \to L_2(\Omega)$ is the 
natural injection throughout this section.

At first we consider the basic example of Section~\ref{Sdtonl2.1}.
The following property of unique continuation plays an important role.

\begin{thm} \label{tdtonl701}
Let $u \in H^1(\Ri^d)$ and $m \in L_\infty(\Ri^d,\Ri)$.
Suppose that $- \Delta u + m \, u = 0$.
If $u$ vanishes on a non-empty open set, then $u = 0$.
\end{thm}
\proof\
See \cite{RS4} Theorem~XIII.57.\hfill$\Box$

\vertspace

Let $m \in L_\infty(\Omega,\Ri)$.
This theorem allows us to prove that $W(\gota) = \{ 0 \} $
for the basic example, where $\gota$ is as in (\ref{exdtonl450;1}).
Recall that  $\gota \colon H^1(\Omega) \times H^1(\Omega) \to \Ci$ is given by
\[
\gota(u,v) 
= \int_\Omega \nabla u \cdot \overline{\nabla v}
   + \int_\Omega m \, u \, \overline v
.  \]

\begin{prop} \label{pdtonl702}
One has $W(\gota) = \{ 0 \} $.
\end{prop}
\proof\
If $f \colon \Omega \to \Ci$ is a function we denote by $\tilde f \colon \Ri^d \to \Ci$ 
the extension of $f$ by $0$.
Note that $\ker j = H^1_0(\Omega)$.
Let $u \in W(\gota)$.
Then $u \in H^1_0(\Omega)$, $- \Delta u + m \, u = 0$ and $\partial_\nu u = 0$.
Since $u \in H^1_0(\Omega)$ one has $\tilde u \in H^1(\Ri^d)$ and 
$\partial_k \tilde u = \widetilde{\partial_k u}$ for all $k \in \{ 1,\ldots,d \} $
cf.\ \cite{Bre} Proposition~IX.18(iii).
Then 
\[
\int_{\Ri^d} \nabla \tilde u \cdot \overline{\nabla v}
   + \int_{\Ri^d} \tilde m \, \tilde u \, \overline v
= \int_\Omega \nabla u \cdot \overline{\nabla v}
   + \int_\Omega m \, u \, \overline v
= 0
\]
for all $v \in \cd(\Ri^d)$.
Therefore $- \Delta \tilde u + \tilde m \, \tilde u = 0$ weakly in $\Ri^d$.
Thus it follows from Theorem~\ref{tdtonl701} that $\tilde u = 0$.\hfill$\Box$

\vertspace

Now we can deduce from Theorem~\ref{tdtonl502.3} the following convergence 
result.

\begin{thm} \label{tdtonl703}
Let $(m_n)_{n \in \Ni}$ be a sequence in $L_\infty(\Omega,\Ri)$ and let 
$m \in L_\infty(\Omega,\Ri)$.
Suppose that $\lim_{n \to \infty} m_n = m$ weak$^*$ in $L_\infty(\Omega,\Ri)$.
Then $\lim_{n \to \infty} (D_{m_n} + i \, s \, I)^{-1} = (D_m + i \, s \, I)^{-1}$
in $\cl(L_2(\Omega))$ for all $s \in \Ri \setminus \{ 0 \} $.
\end{thm}
\proof\
We consider the forms $\gota_n$ and $\gota$ on $H^1(\Omega)$ given by 
\begin{eqnarray}
\gota_n(u,v) 
& = & \int_\Omega \nabla u \cdot \overline{\nabla v}
   + \int_\Omega m_n \, u \, \overline v , \label{etdtonl703;0.8} \\*
\gota(u,v) 
& = & \int_\Omega \nabla u \cdot \overline{\nabla v}
   + \int_\Omega m \, u \, \overline v .
\label{etdtonl703;1}
\end{eqnarray}
Since the sequence $(m_n)$ is uniformly bounded in $L_\infty(\Omega,\Ri)$ 
it follows that the sequence $(\gota_n)$ is uniformly $\tilde j$-elliptic.
Next, let  $v \in V$ and $u,u_1,u_2,\ldots \in H^1(\Omega)$ with $\lim u_n = u$ 
weakly in $H^1(\Omega)$.
Then $\lim u_n = u$ strongly in $L_2(\Omega)$ and 
$\lim u_n \, \overline v = u \, \overline v$ in $L_1(\Omega)$.
Hence $\lim \int_\Omega m_n \, u_n \, \overline v = \int_\Omega m \, u \, \overline v$.
Moreover, $\lim \int_\Omega \nabla u_n \cdot \overline{\nabla v} 
      = \int_\Omega \nabla u \cdot \overline{\nabla v}$.
Therefore $\lim \gota_n(u_n,v) = \gota(u,v)$ and we have shown that $(\gota_n)$ 
converges weakly to $\gota$.
Now Theorem~\ref{tdtonl502.3} gives the result.\hfill$\Box$

\vertspace

Concerning convergence of the semigroups a further condition is needed.
Even if the functions $m_n$ are constant, the graphs $D_{m_n}$ need not be 
uniformly bounded below, so that the semigroups cannot converge strongly 
by the uniform boundedness principle.
We give an example.

\begin{exam} \label{xdtonl704}
Let $(\lambda_n)$ be a strictly increasing sequence in $\Ri$ such that 
$\lim \lambda_n = \lambda_1^D$, where $\lambda^D_1$ is the first eigenvalue
of $-\Delta^D$.
Choose $m_n = - \lambda_n$ constant and 
abbreviate $D_n := D_{m_n}$.
Then $D_n$ is single-valued.
It follows from \cite{ArM2} Propositions~5 and 3 that 
the sequence $(D_n)$ is not uniformly bounded below.
\end{exam}

If, however, the operator $-\Delta^D + m$ is invertible, then the sequence is uniformly 
bounded below. 
This is the content of the next theorem.
Note that $-\Delta^D + m$ is the operator associated with the
classical form $\gota|_{H^1_0(\Omega) \times H^1_0(\Omega)}$ in $L_2(\Omega)$,
where $\gota$ is as in (\ref{etdtonl703;1}).

\begin{thm} \label{tdtonl705}
Let $(m_n)_{n \in \Ni}$ be a sequence in $L_\infty(\Omega,\Ri)$ and let 
$m \in L_\infty(\Omega,\Ri)$.
Suppose that $\lim_{n \to \infty} m_n = m$ weak$^*$ in $L_\infty(\Omega,\Ri)$.
Moreover, assume that $0 \in \rho(-\Delta^D + m)$.
Then the graphs $D_{m_n}$ are bounded below uniformly in $n \in \Ni$ and 
$\lim_{n \to \infty} e^{-t D_{m_n}} = e^{-t D_m}$
in $\cl(L_2(\Gamma))$ for all $t > 0$.
\end{thm}
\proof\
Let $\gota_n$ and $\gota$ be as in (\ref{etdtonl703;0.8}) and (\ref{etdtonl703;1}).
Then $V(\gota) \cap \ker j = \ker(-\Delta^D + m)$.
But $\ker(-\Delta^D + m) = \{ 0 \} $ by assumption.
Therefore $V(\gota) \cap \ker j = \{ 0 \} $ and the graphs $D_{m_n}$ are 
bounded below uniformly in $n \in \Ni$ by Theorem~\ref{tdtonl509.5}.
Now apply Theorems~\ref{tdtonl703} and \ref{tdtonl605}.\hspace*{5mm}\hfill$\Box$

\vertspace

We next consider elliptic operators in divergence form
with real symmetric coefficients instead of the Laplacian.
It is interesting that now the unique continuation property 
depends on the regularity of the coefficients (if $d \geq 3$) and the 
equality $W(\gota) = \{ 0 \} $ is not always valid.

For all $k,l \in \{ 1,\ldots,d \} $ let $a_{kl} \colon \Omega \to \Ri$ 
be bounded, measurable with $a_{kl} = a_{lk}$ and assume that there exists 
a $\mu > 0$ such that the uniform ellipticity condition
\begin{equation}
\sum_{k,l=1}^d a_{kl}(x) \, \xi_k \, \xi_l
\geq \mu \, |\xi|^2
\label{eSdtonl7;2}
\end{equation}
holds for all $x \in \Omega$ and $\xi \in \Ri^d$.
Moreover, let $c \in L_\infty(\Omega,\Ri)$.
We associate with these coefficients the symmetric form 
$\gota \colon H^1(\Omega) \times H^1(\Omega) \to \Ci$ given by 
\begin{equation}
\gota(u,v) 
= \int_\Omega \sum_{k,l=1}^d a_{kl} \, (\partial_k u) \, \overline{\partial_l v}
   + \int_\Omega c \, u \, \overline v
.  
\label{eSdtonl7;3}
\end{equation}
Then $\gota$ is $\tilde j$-elliptic.
We denote by $D_\gota$ the self-adjoint graph associated with $(\gota,j)$.

We next consider the space $W(\gota)$.

\begin{prop} \label{pdtonl706}
Assume that the second-order coefficients $a_{kl}$ are Lipschitz
continuous for all $k,l \in \{ 1,\ldots,d \} $ or that $d = 2$.
Then $W(\gota) = \{ 0 \} $.
\end{prop}
\proof\
By the McShane--Whitney extension theorem each Lipschitz continuous 
function on $\Omega$ has a Lipschitz continuous extension to $\Ri^d$.
Hence we may assume that the $a_{kl}$ are defined on $\Ri^d$, 
are Lipschitz continuous, symmetric, bounded and satisfy the ellipticity
condition (\ref{eSdtonl7;2}) uniformly for all $x \in \Ri^d$, possibly with a 
different value of $\mu > 0$.
Moreover, we may assume that $a_{kl}(x) = \delta_{kl}$ for all 
$x \in \Ri^d$ with $|x|$ large.
Arguing as in the proof of Proposition~\ref{pdtonl702}
and using the unique continuation property
(see \cite{AKS} Section 5, Remark 3, or (W) in the Introduction of \cite{Kur},
or \cite{GL} Theorem~1.1)
one deduces that $W(\gota) = \{ 0 \} $.

In case $d=2$ Schulz \cite{Schulz} proved the unique continuation property 
without the assumption that the leading coefficients are Lipschitz 
continuous.\hfill$\Box$

\vertspace

Triviality of $W(\gota)$ is a most interesting property.
In fact, given $s \in \Ri \setminus \{ 0 \} $ and $h \in L_2(\Gamma)$ we 
can always find a $u \in H^1(\Omega)$ such that 
\[
\gota(u,v) + i \, s \int_\Gamma \Tr u \, \overline{\Tr v} 
= \int_\Gamma h \, \overline{\Tr v} 
\]
for all $v \in H^1(\Omega)$.
By our results the trace of $u$ does not depend on the choice of $u$
and by definition
$\Tr u = ( D_\gota + i \, s \, I)^{-1} h$.
The element $u$ is unique if and only if $W(\gota) = \{ 0 \} $
by (\ref{etdton202;1}).

Filonov \cite{Fil1} constructed a remarkable example of an elliptic operator 
with H\"older continuous real symmetric coefficients on the open ball
$\Omega$ in $\Ri^3$
(even H\"older continuous of order $\nu$ for every $\nu \in (0,1)$)
and a $w \in \cd(\Omega)$ such that $w \in W(\gota) \setminus \{ 0 \} $.
Thus we have non-uniqueness of the function $u$ above.
By the theorem on unique continuation, these coefficients cannot be
Lipschitz continuous.

For Lipschitz continuous coefficients in the limit we can prove the 
following convergence result.

\begin{thm} \label{tdtonl707}
For all $k,l \in \{ 1,\ldots,d \} $ and $n \in \Ni$ let 
$a^{(n)}_{kl} \in L_\infty(\Omega,\Ri)$.
Further assume $d=2$ and $a_{kl} \in L_\infty(\Omega,\Ri)$, or
$d \geq 3$ and $a_{kl} \in W^{1,\infty}(\Omega,\Ri)$
for all $k,l \in \{ 1,\ldots,d \} $.
Next let $c_n,c \in L_\infty(\Omega,\Ri)$
and fix $\mu > 0$.
Suppose that 
$a^{(n)}_{kl} = a^{(n)}_{lk}$ for all $k,l \in \{ 1,\ldots,d \} $
and 
\[
\sum_{k,l=1}^d a^{(n)}_{kl}(x) \, \xi_k \, \xi_l
\geq \mu \, |\xi|^2
\]
for all $x \in \Omega$, $\xi \in \Ri^d$ and $n \in \Ni$.
Suppose that $\lim_{n \to \infty} \|a^{(n)}_{kl} - a_{kl}\|_\infty = 0$ 
for all $k,l \in \{ 1,\ldots,d \} $
and $\lim_{n \to \infty} c_n = c$ weak$^*$ in $L_\infty(\Omega)$.
For all $n \in \Ni$ define the form 
$\gota_n \colon H^1(\Omega) \times H^1(\Omega) \to \Ci$
by 
\[
\gota_n(u,v) 
= \int_\Omega \sum_{k,l=1}^d a^{(n)}_{kl} \, (\partial_k u) \, \overline{\partial_l v}
   + \int_\Omega c_n \, u \, \overline v
\]
and define $\gota$ as in {\rm (\ref{eSdtonl7;3})}.
Let $A_n$ be the graph associated with $(\gota_n,j)$ and $A$ the graph 
associated with $(\gota,j)$.
Then 
\[
\lim_{n \to \infty} (A_n + i \, s \, I)^{-1} = (A + i \, s \, I)^{-1}
\]
in $\cl(L_2(\Gamma))$ for all $s \in \Ri \setminus \{ 0 \} $.
\end{thm}
\proof\
By Proposition~\ref{pdtonl706} we know that $W(\gota) = \{ 0 \} $.
Arguing as in the proof of Theorem~\ref{tdtonl703} one deduces that
the sequence $(\gota_n)_{n \in \Ni}$ converges weakly to the form $\gota$.
Now the claim is a consequence of Theorem~\ref{tdtonl502.3}.\hfill$\Box$

\vertspace

Adopt the assumptions and notation of Theorem~\ref{tdtonl707}.
Let $A^D$ be the operator in $L_2(\Omega)$ associated with the form
$\gota|_{H^1_0(\Omega) \times H^1_0(\Omega)}$. 
Assuming that $A^D$ is invertible we can deduce uniform convergence 
of the associated semigroups exactly as in Theorem~\ref{tdtonl705}.

\smallskip

It is remarkable that we can deduce from Proposition~\ref{pdtonl506}
that the set of all second-order coefficients for which $W(\gota) = \{ 0 \} $
is open in the following sense.
Let $\mu > 0$.
Consider the set 
\[
Q = \{ (c_{kl})_{kl} \in L_\infty(\Omega,\Ri)^{d \times d}_{\rm sym}
   : \sum_{k,l=1}^d c_{kl}(x) \, \xi_k \, \xi_l \geq \mu \, |\xi|^2 
      \mbox{ for all } \xi \in \Ri^d \mbox{ and } x \in \Omega \}
.  \]
To each $c = (c_{kl})_{kl} \in Q$ we associate the form 
$\gota_c \colon H^1(\Omega) \times H^1(\Omega) \to \Ci$ defined by
\[
\gota_c(u,v) = \sum_{k,l=1}^d \int_\Omega c_{kl} \, (\partial_k u) \, \overline{\partial_l v}
.  \]
Let $c \in Q$ and suppose that $W(\gota_c) = \{ 0 \} $.
Then there exists an $\varepsilon > 0$ such that 
for all $\hat c \in Q$ with $\|c - \hat c\|_{L_\infty(\Omega)^{d \times d}} < \varepsilon$
one also has the uniqueness property $W(\gota_{\hat c}) = \{ 0 \} $.

\section{Lumer--Phillips by hidden compactness} \label{Sdtonl4}

In this section we replace the condition that the form $\gota$ is 
symmetric with the condition that $\gota$ is accretive. 

We say that a graph $A$ is {\bf accretive} if $\RRe (x,y) \geq 0$ for all 
$(x,y) \in A$.
The graph $A$ is called {\bf $m$-accretive} if it is accretive and 
$A + I$ is surjective.
Our point here is that this latter condition of surjectivity can be 
obtained by hidden compactness.
More precisely, we show the following.

\begin{thm} \label{tdtonl430}
Let $V$ and $H$ be Hilbert spaces.
Let $\gota \colon V \times V \to \Ci$ be an accretive continuous sesquilinear form.
Further let $j \in \cl(V,H)$.
Let $A$ be the graph associated with 
$(\gota,j)$.
If $\gota$ is compactly elliptic, then $A$ is $m$-accretive.
\end{thm}
\proof\
Clearly the graph $A$ is accretive.

First suppose that $W(\gota) = \{ 0 \} $.
We shall show that $I + A$ is surjective.
Define the sesquilinear form $\gotb \colon V \times V \to \Ci$ by 
\[
\gotb(u,v) = \gota(u,v) + (j(u), j(v))_H
.  \]
Then $\gotb$ is  compactly elliptic.
Define $\cb \colon V \to V'$ by 
$(\cb u,v)_{V' \times V} = \gotb(u,v)$.
We shall show that $\cb$ is injective.
Let $u \in V$ and suppose that $\cb u = 0$.
Then 
\[
0 = \RRe (\cb u, u)
= \RRe \gota(u) + \|j(u)\|_H^2
\geq \|j(u)\|_H^2
.  \]
So $j(u) = 0$.
Then for all $v \in V$ one has 
\[
0 
= (\cb u,v)_{V' \times V}
= \gotb(u,v)
= \gota(u,v) + (j(u), j(v))_H
= \gota(u,v)
.  \]
So $u \in W(\gota) = \{ 0 \} $ by assumption.
Therefore $\cb$ is injective and hence also surjective by the 
Fredholm--Lax--Milgram lemma, Lemma~\ref{ldtonl223}.
Now let $y \in H$.
Define $\alpha \colon V \to \Ci$ by 
$\alpha(v) = (y, j(v))_H$.
Then $\alpha \in V'$ since $j$ is continuous.
Because $\cb$ is surjective, there exists a (unique) $u \in V$ such that
$\cb u = \alpha$.
Then for all $v \in V$ one has
\begin{eqnarray*}
(y, j(v))_H
& = & (\cb u,v)_{V' \times V}
= \gotb(u,v)
= \gota(u,v) + (j(u), j(v))_H
= \gota(u,v) + (x, j(v))_H
,
\end{eqnarray*}
where $x = j(u)$.
So $x \in D(A)$ and $(A + I)x = y$.
This proves that $A$ is $m$-accretive if $W(\gota) = \{ 0 \} $.

Finally we drop the assumption that $W(\gota) = \{ 0 \} $.
There exists a unique $T \in \cl(V)$ such that 
$\gota(u,v) = (Tu,v)_V$ for all $u,v \in V$.
Then $T$ is $m$-accretive.
So $\ker T = \ker T^*$.
Hence $W(\gota) = \ker j \cap \ker T = \ker j \cap \ker T^* = W(\gota^*)$.
Let $V_1 = W(\gota)^\perp$, where the orthogonal complement is in $V$.
Define $\gota_1 = \gota|_{V_1 \times V_1}$ and $j_1 = j|_{V_1}$.
Then $\gota_1$ is compactly elliptic too by Proposition~\ref{pdtonl223.7}\ref{pdtonl223.7-3}.
Let $u \in W(\gota_1)$.
Then $u \in V_1$, $j(u) = 0$ and $\gota(u,v) = 0$ for all $v \in V_1$.
If $w \in W(\gota)$ then $w \in W(\gota^*)$ and $\gota^*(w,u) = 0$ by definition of $W(\gota^*)$.
So $\gota(u,w) = \overline{\gota^*(w,u)} = 0$.
Hence by linearity $\gota(u,v) = 0$ for all $v \in V$.
Therefore $u \in W(\gota)$.
So $u \in W(\gota) \cap V_1 \subset W(\gota) \cap W(\gota)^\perp = \{ 0 \} $.
Thus $W(\gota_1) = \{ 0 \} $.
Let $A_1$ be the operator associated to $(\gota_1,j_1)$.
By the first part of the proof, the operator $A_1$ is $m$-accretive.
Since $A$ is accretive, it suffices to show that $A_1 \subset A$.
Let $(x,y) \in A_1$.
By definition there exists a $u \in V_1$ such that $j_1(u) = x$ and 
$\gota_1(u,v) = (y,j_1(v))_H$ for all $v \in V_1$.
Let $w \in W(\gota)$.
Then $w \in W(\gota^*)$, so as above one deduces that $\gota(u,w) = 0$.
So $\gota(u,w) = (y,0)_H = (y,j(w))_H$.
Then by linearity one has $\gota(u,v) = (y,j(v))_H$ for all $v \in V$.
Therefore $(x,y) \in A$ and $A_1 \subset A$.
This completes the proof of the theorem.\hfill$\Box$

\begin{exam} \label{xdtonl432}
Let $V$ and $H$ be Hilbert spaces such that $V$ is densely and compactly 
embedded in $H$.
Let $\gota \colon V \times V \to \Ci$ be a continuous accretive form which 
is $H$-elliptic.
Let $A$ be the $m$-accretive {\em operator} associated with $\gota$.
Next, let $\widehat H$ be a Hilbert space and $T \in \cl(H,\widehat H)$.
Define the graph $B$ in $\widehat H$ by 
\[
B = \{ (x,y) \in \widehat H \times \widehat H: \mbox{there exists a } u \in D(A)
      \mbox{ such that } T u = x \mbox{ and } T^* y = A u \}
.  \]
Then $B$ is $m$-accretive.
The proof is as follows.
Define $j \colon V \to \widehat H$ by $j(u) = Tu$.
Let $B_2$ be the graph associated with $(\gota,j)$.
Then $B_2$ is $m$-accretive by Theorem~\ref{tdtonl430}.

We show that $B$ is accretive.
Let $(x,y) \in B$.
Then there exists a $u \in D(A)$ such that $T u = x$ and 
$Au = T^* y$.
Therefore 
\[
\RRe (x,y)_{\widehat H}
= \RRe (T u,y)_{\widehat H}
= \RRe (u,T^* y)_H
= \RRe (u,Au)
= \RRe \overline{\gota(u)}
\geq 0
.  \]
So $B$ is accretive.

Next we show that $B_2 \subset B$.
Let $(x,y) \in B_2$.
Then there exists a $u \in V$ such that $j(u) = x$ and
$\gota(u,v) = (y,j(v))_{\widehat H}$ for all $v \in V$.
Hence $\gota(u,v) = (y,T v)_{\widehat H} = (T^* y, v)_H$ for all $v \in V$.
This implies that $u \in D(A)$ and $Au = T^* y$.
Therefore $(x,y) \in B$.

Since an $m$-accretive operator does not have a proper accretive 
extension, it follows that $B = B_2$.
\end{exam}

In Theorem~\ref{tdtonl430} the condition that $\gota$ is accretive 
cannot be removed in general. 
An example is as follows.

\begin{exam} \label{xdtonl435}
Let $V = \widetilde H = \Ci$ and $H = \Ci$.
Let $\gota \colon V \times V \to \Ci$ be given by $\gota(u,v) = u_2 \, \overline{v_1}$.
Define $j(u) = u_1$ and $\tilde j(u) = u$.
Then $\gota$ is $\tilde j$-elliptic.
The graph associated with $(\gota,j)$ is $\Ci \times \Ci$,
which is not lower-bounded in the sense that there is no $M > 0$ such that 
$\RRe (x,y)_H \geq -M \, \|x\|_H^2$ for all $(x,y) \in A$.

Note that for symmetric $\gota$ we proved lower boundedness in Theorem~\ref{tdtonl210.5}.
\end{exam}

We conclude with an example of a form $\gota$ and a map $j$ for which one has 
hidden compactness, but such that the form $\gota$ is not $j$-elliptic.

\begin{exam} \label{xdtonl434}
Let $\Omega$ be a Lipschitz domain such that the boundary $\Gamma$ has measure~$1$.
Choose $V = H^1(\Omega)$, $H = L_2(\Gamma)$ and 
$\gota(u,v) = \int_\Omega \nabla u \cdot \overline{\nabla v}$.
Moreover, define $B \colon L_2(\Gamma) \to L_2(\Gamma)$ by 
$B g = g - (g,\one_\Gamma)_{L_2(\Gamma)} \, \one_\Gamma$.
Then $B \one_\Gamma = 0$.
Choose $j = B \circ \Tr$.
Then $\gota$ is accretive and symmetric.
Let $A$ be the graph associated with $(\gota,j)$.
Then $A$ is self-adjoint and $m$-accretive by hidden compactness.
But the form $\gota$ is not $j$-elliptic since
$\gota(\one_\Omega) + \omega \, \|j(\one_\Omega)\|_{L_2(\Gamma)}^2 = 0$
for all $\omega \in \Ri$.
\end{exam}

\subsection*{Acknowledgements}
The first-named author is most grateful for a most stimulating stay at the 
University of Auckland and the great hospitality in Auckland.
The second-named author is most grateful for the hospitality extended
to him during a fruitful stay at the University of Ulm.
He wishes to thank the University of Ulm for financial support.
The third-named author is supported by a fellowship of the
Alexander von Humboldt Foundation, Germany.
Part of this work is supported by the Marsden Fund Council from Government funding,
administered by the Royal Society of New Zealand.

\end{document}